\documentclass[11pt]{article}
\begin{document}
\newtheorem{defn0}{Definition}[section]
\newtheorem{prop0}[defn0]{Proposition}
\newtheorem{thm0}[defn0]{Theorem}
\newtheorem{lemma0}[defn0]{Lemma}
\newtheorem{coro0}[defn0]{Corollary}
\newtheorem{exa}[defn0]{Example}
\newtheorem{rem0}[defn0]{Remark}
\def\rig#1{\smash{ \mathop{\longrightarrow}
    \limits^{#1}}}
\def\nwar#1{\nwarrow
   \rlap{$\vcenter{\hbox{$\scriptstyle#1$}}$}}
\def\near#1{\nearrow
   \rlap{$\vcenter{\hbox{$\scriptstyle#1$}}$}}
\def\sear#1{\searrow
   \rlap{$\vcenter{\hbox{$\scriptstyle#1$}}$}}
\def\swar#1{\swarrow
   \rlap{$\vcenter{\hbox{$\scriptstyle#1$}}$}}
\def\dow#1{\Big\downarrow
   \rlap{$\vcenter{\hbox{$\scriptstyle#1$}}$}}
\def\up#1{\Big\uparrow
   \rlap{$\vcenter{\hbox{$\scriptstyle#1$}}$}}
\def\lef#1{\smash{ \mathop{\longleftarrow}
    \limits^{#1}}}
\def\O{{\cal O}}
\def\L{{\cal L}}
\def\P#1{{\bf P}^#1}
\newcommand{\defref}[1]{Definition~\ref{#1}}
\newcommand{\propref}[1]{Proposition~\ref{#1}}
\newcommand{\thmref}[1]{Theorem~\ref{#1}}
\newcommand{\lemref}[1]{Lemma~\ref{#1}}
\newcommand{\corref}[1]{Corollary~\ref{#1}}
\newcommand{\exref}[1]{Example~\ref{#1}}
\newcommand{\secref}[1]{Section~\ref{#1}}

\newcommand{\qedd}{\hfill\framebox[2mm]{\ }\medskip}
\newcommand{\codim}{\textrm{codim}}
\author{Maria Chiara Brambilla and Giorgio Ottaviani}
\title{On the Alexander-Hirschowitz Theorem}
\date{}
\maketitle
\abstract{
The Alexander-Hirschowitz theorem says that
a general collection of $k$ double points in ${\bf P}^n$ imposes
independent conditions on homogeneous polynomials of degree $d$
with a well known list of exceptions.
Alexander and Hirschowitz completed its proof in 1995, solving a long standing classical problem,
connected with the Waring problem for polynomials.
We expose a self-contained  proof based mainly on previous works by
Terracini,  Hirschowitz, Alexander and Chandler, with a few
simplifications. We claim originality only in the case $d=3$,
where our proof is shorter. We end with an account of the history
of the work on this problem.
}
\medskip

\noindent{\it AMS Subject Classification:} 01-02, 14C20, 15A72, 14M17\\
Both authors are partially supported by Italian MUR and are members of GNSAGA.

\section{Introduction}
The aim of this paper is to expose a proof of the following theorem.

\begin{thm0}[Alexander-Hirschowitz]\label{AH}
Let $X$ be a general collection of  $k$ double points in ${\bf P}^n={\bf
P}(V)$ (over an algebraically closed field of characteristic zero) and let $S^dV^{\vee}$
be the space of homogeneous polynomials of degree $d$. Let
$I_X(d)\subseteq S^dV^{\vee}$ be the subspace of polynomials
through $X$, that is with all first partial derivatives vanishing
at the points of $X$. Then the subspace $I_X(d)$ has the expected
codimension $\min\left((n+1)k,{{n+d}\choose{n}}\right)$ except in
the following cases

 $\bullet\quad d=2,2\leq k\leq n$;

$\bullet\quad n=2,d=4,k=5;$

 $\bullet\quad n=3,d=4,k=9;$

 $\bullet\quad n=4,d=3,k=7;$

 $\bullet\quad n=4,d=4,k=14.$
\end{thm0}

We remark that the case $n=1$ is the only one where the assumption that $X$
is general is not necessary.

More information on
the exceptional cases is contained in Section \ref{sezione3}.

This theorem has an equivalent formulation in terms of higher secant
varieties. Given a projective variety $Y$, the $k$-secant variety
$\sigma_k(Y)$ is the Zariski closure of the union of all the
linear spans $< p_1,\ldots ,p_k >$ where $p_i\in Y$ (see \cite{Ru} or \cite{Z}). In
particular $\sigma_1(Y)$ coincides with  $Y$ and $\sigma_2(Y)$ is the usual secant
variety. Consider the Veronese embedding $V^{d,n}\subset{\bf P}^m$
of degree $d$ of ${\bf P}^n$, that is the image of the linear
system given by all homogeneous polynomials of degree $d$, where
$m={{n+d}\choose{n}}-1$. It is easy to check that
$\dim\sigma_k(V^{d,n})\le\min\left((n+1)k-1,m\right)$ and when
the equality holds we say that $\sigma_k(V^{d,n})$ has the
expected dimension.

\begin{thm0}[Equivalent formulation of Theorem \ref{AH}]\label{AH2}
The higher secant variety $\sigma_k(V^{d,n})$ has the expected dimension
with the same exceptions of \thmref{AH}.
\end{thm0}

\thmref{AH2} still holds if the characteristic of the base field ${\bf K}$ is
bigger than $d$ and $d>2$ (\cite[Corollary I.62]{IK}),
but the case  $\textrm{char} ({\bf K})=d$ is open as far as we know.
The equivalence between \thmref{AH} and \thmref{AH2} holds if  $\textrm{char}( {\bf K})=0$, and since we want to switch
freely between the two formulations we work with this assumption. Let us mention
that in \cite{AH1} \thmref{AH} is stated with the weaker assumption that ${\bf K}$ is infinite.

Since the general element in $\sigma_k(V^{d,n})$ can be expressed
as the sum of $k$ $d$-th powers of linear forms, a consequence of \thmref{AH2}
is that the general homogeneous polynomial of degree $d$ in $n+1$ variables can be expressed as the sum
of $\lceil\frac{1}{n+1}{{n+d}\choose{d}}\rceil$ $d$-th powers of linear forms
with the same list of exceptions (this is called the Waring problem
for polynomials, see \cite{IK}).

In the case $n=1$, the Veronese embedding $V^{d,1}$ is the rational normal curve and
there are no exceptions at all. The case $n=2$ was proved  by
Campbell \cite{Ca},  Palatini \cite{Pa2} and Terracini
\cite{Te2}, see the historical section \ref{sezione7}. In
\cite{Pa2} Palatini stated Theorem \ref{AH} as a plausible
conjecture. In \cite{Te1} Terracini proved his famous two
``lemmas'', which turned out to be crucial keys to solve the
general problem. In 1931 Bronowski claimed to have a proof of
Theorem \ref{AH}, but his proof was fallacious. Finally the proof
was found in 1995 by Alexander and Hirschowitz along a series of
brilliant papers, culminating with \cite{AH1}  so
that Theorem \ref{AH} is now called the Alexander-Hirschowitz
theorem. They introduced the so called differential Horace's method to
attack the problem. The proof was simplified in \cite{AH2}.
 In 2001 K. Chandler achieved a further
simplification in \cite{Ch} and \cite{Ch2}.   The higher multiplicity case is still open and it is  a subject of active research,
due to a striking conjecture named after Segre-Gimigliano-Harbourne-Hirschowitz,
see \cite{Ci} for a survey.

In 2006 we ran a seminar in Firenze trying to understand this problem. This note is a result of that
seminar, and reflects the historical path that we have chosen.
We are able to present a self-contained and detailed proof of the Alexander-Hirschowitz theorem, starting from scratch, with several simplifications
on the road tracked by Terracini, Hirschowitz, Alexander and Chandler.

The reader already accustomed to this topic can skip
Section \ref{sezione4} which is added only to clarify the problem
and jump directly
to Sections \ref{sezione5} and \ref{horace}, which contain our original contributions
(especially Section \ref{sezione5} about cubics,
while in Section \ref{horace} we supplied \cite{Ch} with more details).

The Veronese varieties are one of the few classes
of varieties where the dimension of the higher secant varieties is completely known.
See \cite{CGG1}, \cite{CGG2}, \cite{MG}, \cite{AOP} for related work on Segre and Grassmann varieties.

We thank all the participants to the seminar for their criticism, especially Luca Chiantini. We thank also Ciro Ciliberto
for his remarks concerning the historical Section \ref{sezione7} and
 Edoardo Ballico for helpful comments on the previous version of our paper. 

\section{Notation and Terracini's two lemmas}
\label{notation}
For any real number $x$, $\lfloor x \rfloor$ is the greatest integer smaller than or
equal to $x$,
$\lceil x \rceil$  is the smallest integer greater than or equal to $x$.
Let $V$ be a vector space of dimension $n+1$ over an algebraically closed field ${\bf K}$
of characteristic zero.
 Let
${\bf P}^n={\bf P}(V)$ be the projective space of lines in $V$.
If $f\in V\setminus\{0\}$ we denote by $[f]$ the line spanned by $f$
and also the corresponding point in ${\bf P}(V)$. Let $S=\oplus_d S^dV$
be the symmetric algebra of $V$ and $S^{\vee}=\oplus_d S^dV^{\vee}$ its dual.
We have the natural pairing $S^dV\otimes S^dV^{\vee}\to {\bf K}$ which we denote by $(\ ,\ )$.
Then $S^dV^{\vee}$ is the space of homogeneous polynomials over
${\bf P}(V)$ and a polynomial $h\in S^dV^{\vee}$ vanishes at $[f]\in {\bf P}(V)$
if and only if  $(f^d,h)=0$. The {\it Veronese variety} $V^{d,n}$
is the image of the embedding $[v]\mapsto [v^d]$
of ${\bf P}(V)$ in ${\bf P}(S^dV)={\bf P}^m$, where $m={{n+d}\choose{n}}-1$.
If $f\in V$, it is easy to check that the projective tangent space $T_{[f^d]}V^{d,n}\subseteq {\bf P}(S^dV)$
is equal to $\{[f^{d-1}g]| g\in V\}$ (to see this, compute the Taylor expansion of $(f+\epsilon g)^d$ at $\epsilon=0$).

The maximal ideal corresponding to $f\in V$ is
$${ m}_{[f]}:=\{h\in  S^{\vee}|h(f)=0\}.$$
It contains all the hypersurfaces which pass through $[f]$.
Its power ${ m}^2_{[f]}$ contains all the hypersurfaces which are singular at $[f]$,
it defines a scheme which is denoted as $[f]^2$ and it is called a {\it double point}.
Note that a hypersurface is singular at $[f]$ if and only if it contains $[f]^2$.

In order to state the relation between the higher secant varieties to the Veronese varieties and the double points of hypersurfaces
we need the following proposition, well known to Palatini and Terracini, usually attributed to Lasker
\cite{La}, the Hilbert's student who proved the primary decomposition for ideals in polynomial rings and is widely known as
chess world champion at the beginning of XX century.

\begin{prop0} \label{lasker}{\bf (Lasker)}
Given $T_{[f^d]}V^{d,n}\subseteq {\bf P}(S^dV)$, its (projectivized) orthogonal
$\left(T_{[f^d]}V^{d,n}\right)^{\perp}\subseteq {\bf P}(S^dV^{\vee})$ consists of
all the hypersurfaces singular at $[f]$.
More precisely, if we denote by $C(V^{d,n})$ the affine cone over $V^{d,n}$, then the following holds
$$\left(T_{f^d}C(V^{d,n})\right)^{\perp} = \left({ m}^2_{[f]}\right)_d\subseteq S^dV^{\vee}$$
\end{prop0}
{\it Proof.} Let $e_0,\ldots e_n$ be a basis of $V$ and
$x_0,\ldots ,x_n$ its dual basis. Due to the $GL(V)$-action it is enough to check the statement for
$f=e_0$. Then  ${ m}_{[f]}=(x_1,\ldots ,x_n)$, ${ m}^2_{[f]}=(x_1^2,x_1x_2,\ldots ,x_n^2)$,
so that $\left({ m}^2_{[f]}\right)_d$  is generated by all monomials of degree $d$ with the exception of
$x_0^d,x_0^{d-1}x_1,\ldots ,x_0^{d-1}x_n$.\\ Since  $T_{e_0^d}C(V^{d,n})=<e_0^d, e_0^{d-1}e_1,\ldots ,e_0^{d-1}e_n>$
the thesis follows.\qedd

\begin{lemma0} {\bf (First Terracini lemma)} \label{ter1}
Let  $p_1,\ldots p_k\in Y$ be general points and  \linebreak
 $z\in <p_1,\ldots ,p_k>$
 a general point. Then
$$T_z\sigma_k(\,Y)=<T_{p_1}Y,\ldots ,T_{p_k}Y>.$$
\end{lemma0}
{\it Proof.} Let $Y(\tau)=Y(\tau_1,\ldots ,\tau_n)$ be a local parametrization of $Y$.
We denote by $Y_j(\tau )$ the partial derivative with respect to $\tau_j$.
Let $p_i$ be the point corresponding to $\tau^i=(\tau^i_1,\ldots,\tau^i_n )$.
The space $<T_{p_1}Y,\ldots ,T_{p_k}Y>$ is spanned by the $k(n+1)$ rows
of the following matrix
$$\begin{array}{c}
\vdots\\
Y(\tau^i)\\
Y_1(\tau^i)\\
\vdots     \\
Y_n(\tau^i)\\
\vdots\\
\end{array}$$
(here we write only the $i$-th block of rows, $i=1,\ldots ,k$).

We write also the local parametrization of  $\sigma_k(Y)$ given by
$$\Phi(\tau^1,\ldots ,\tau^k,\lambda_1,\ldots ,\lambda_{k-1})=
\sum_{i=1}^{k-1}\lambda_iY(\tau^i)+Y(\tau^k)$$
depending on $kn$ parameters $\tau^i_j$ and $k-1$ parameters $\lambda_i$.
The matrix whose rows are given by $\Phi$ and its $kn+k-1$ partial
derivatives computed at $z$ is
$$\begin{array}{c}
\sum_{i=1}^{k-1}\lambda_iY(\tau^i)+Y(\tau^k)\\
\vdots\\
\lambda_iY_1(\tau^i)\\
\vdots     \\
\lambda_iY_n(\tau^i)\\
\vdots\\
Y_1(\tau^k)\\
\vdots     \\
Y_n(\tau^k)\\
Y(\tau^1)\\
\vdots\\
Y(\tau^{k-1})\\
\end{array}$$
and its rows span $T_z\sigma_k(Y)$.
It is elementary to check that the two above matrices are obtained
one from the other by performing elementary operations on rows, hence they have the same row space and the same
rank. \qedd

The \propref{lasker} and \lemref{ter1} allow to prove the equivalence between
\thmref{AH} and \thmref{AH2}. Indeed let $X=\{ p_1^2,\ldots ,p_k^2\}$ be a collection
of double points in $\P n$ and choose some representatives $v_i\in V$ such that $[v_i]=p_i$ for
$i=1,\ldots ,k$. The subspace
$$I_X(d)=\bigcap_{i=1}^k\left[{ m}^2_{p_i}\right]_d$$
is equal by \propref{lasker} to
$$\bigcap_{i=1}^k \left(T_{v_i^d}C(V^{d,n})\right)^{\perp}=\left( <T_{v_1^d}C(V^{d,n}),\ldots , T_{v_k^d}C(V^{d,n}) >\right)^{\perp}\subseteq
S^dV^{\vee}$$
so that its codimension is equal to the dimension of
$$<T_{v_1^d}C(V^{d,n}),\ldots , T_{v_k^d}C(V^{d,n})>\subseteq S^dV$$
which in turn is equal to
$$\dim <T_{[v_1^d]}V^{d,n},\ldots , T_{[v_k^d]}V^{d,n} >+1$$
where we consider now the projective dimension.
Summing up, by using \lemref{ter1}, the genericity assumption on the points
and the fact that $\sigma_k (V^{d,n})$ is an irreducible variety, we get
$$\textrm{codim }I_X(d)=\textrm{dim } \sigma_k (V^{d,n})+1$$
and the equivalence between \thmref{AH} and \thmref{AH2} is evident from this equality.
 \medskip

We say that a collection $X$ of double points {\em imposes independent conditions
on ${\cal O}_{\P n}(d)$} if
the codimension of $I_X(d)$ in $S^dV^{\vee}$ is $\min\{ {{n+d}\choose{n}}, k(n+1)\}$.
It always holds $\textrm{codim}\ I_X(d)\le \min\{ {{n+d}\choose{n}}, k(n+1)\}$.
Moreover if  $\textrm{codim}\ I_X(d)=  k(n+1)$ and $X'\subset X$ is a collection
of $k'$ double points then $\textrm{codim}\ I_{X'}(d)=  k'(n+1)$.
On the other hand if  $\textrm{codim}\ I_X(d)=  {{n+d}\choose{n}}$ and $X''\supset X$ is a collection
of $k''$ double points then $\textrm{codim}\ I_{X''}(d)=  {{n+d}\choose{n}}$.

\begin{lemma0} {\bf (Second Terracini lemma)}   \label{ter2}
Let $X$ be a  union of double points supported on $p_i$, $i=1,\ldots ,k$. 
We identify the points $p_i$ with their images
on $V^{d,n}$ according to the Veronese embedding. Assume that $X$ does not impose independent conditions
on hypersurfaces of degree $d$.   Then there is a positive dimensional variety $C\subseteq V^{d,n}$
 through $p_1,\ldots p_k$
such that if $p\in C$ then
$T_pV^{d,n}\subseteq <T_{p_1}V^{d,n},\ldots , T_{p_k}V^{d,n}>$.
In particular, by \propref{lasker}, every hypersurface of degree $d$ which is singular at $p_i$ is also singular along $C$.
\end{lemma0}

{\it Proof.}  Let $z$ be a general point in $<p_1,\ldots ,p_k>$.  By Lemma \ref{ter1} we have
$$T_z\sigma_k(V^{d,n})=<T_{p_1}V^{d,n},\ldots ,T_{p_k}V^{d,n}>$$

The secant variety $\sigma_k(V^{d,n})$ is obtained by projecting on the last factor
the abstract secant variety
$\sigma^k(V^{d,n})\subseteq V^{d,n}\times\ldots\times V^{d,n}\times\P m $
which is defined as follows
$$\sigma^k(V^{d,n}):=\overline{\{(q_1,\ldots ,q_k,z) |
z\in <q_1,\ldots ,q_k>, \dim <q_1,\ldots ,q_k>=k-1\} }$$
and has dimension $nk+(k-1)$.

By assumption the dimension of $\sigma_k(V^{d,n})$ is smaller than expected.
Then the fibers $Q_z$ of the above projection have positive dimension and are invariant  under permutations of the first
$k$ factors. Note that $(p_1,\ldots ,p_k)\in Q_z$ and moreover \linebreak
$z\in <q_1,\ldots ,q_k>$ for all $(q_1,\ldots ,q_k)\in Q_z$ such that $\dim <q_1,\ldots ,q_k>=k-1$.
In particular for any such $q_1$ we have that
$T_{q_1}V^{d,n}\subseteq  <T_{p_1}V^{d,n},\ldots , T_{p_k}V^{d,n}>$.

The image of $Q_z$ on the first (or any) component is the variety $C$ we looked for.
 \qedd

{\bf Remark.} It should be mentioned that Terracini proved also a bound on the linear span
of $C$, for details see \cite{CC}. The proofs of the two Lemmas that we have exposed are taken from \cite{Te1}.

\medskip

The first application given by Terracini is the following version of Theorem \ref{AH} in the case $n=2$
(see also the historical section \ref{sezione7}).
\begin{thm0}\label{pidue}
A general union of double points $X\subseteq{\bf P}^2$
imposes independent conditions on plane curves of degree $d$ with the only two exceptions

$d=2,\quad X$\ given by two double points;

$d=4,\quad X$\ given by five double points.
\end{thm0}
{\it Proof.}  We first check the statement for small values of $d$.
It is elementary for $d\le 2$.
Now, every cubic with two double points contains the line
through these two points (by B\'ezout theorem), hence every cubic with three double points is the union of three lines.
It follows easily that the statement is true for $d=3$.
For $d=4$ remind that any quartic with four double points contains a conic through
these points (indeed impose to the conic to pass through a further point and apply
B\'ezout theorem).
Hence there is a unique quartic through five double points,
which is the double conic.

Assume that a general union $X$ of $k$ double points  does not impose independent conditions on plane
 curves of degree $d$.
If $F$ is a plane curve of degree $d$ through $X$, then by Lemma \ref{ter2}
$F$ contains a double curve of degree $2l$ through $X$.
Hence we have the inequalities
$$2l\le d\qquad \textrm{ and  }\qquad k\le\frac{l(l+3)}{2}.$$
We may also assume $$\left\lfloor\frac{1}{3}{{d+2}\choose 2}\right\rfloor\le k$$
because the left-hand side is the maximum expected number of double points imposing independent conditions
on plane curves of degree $d$,
so that we get the inequality
$$\left\lfloor\frac{(d+2)(d+1)}{6}\right\rfloor\le \frac{d}{4}({\frac{d}{2}}+3)$$
which gives $d\le 4$ (already considered) or $d=6$. So the theorem is proved for any $d\neq 6$.
In the case $d=6$ the last inequality is an equality which forces $k=9$.
It remains to prove that the unique sextic which is singular at $9$ general points is the double cubic through
 these points, which follows
again by Lemma \ref{ter2}.
 \qedd

\section{The exceptional cases}
\label{sezione3}
Two double points do not impose independent conditions to the linear system of
quadrics.
Indeed the system of quadrics singular at two points consists of cones having
the vertex containing the line joining the two points,
which has projective dimension ${n\choose 2}>{{n+2}\choose 2}-2(n+1)$.
The same argument works for $k$ general points, $2\leq k\leq n$. In the border
case $k=n$,
the only surviving quadric is the double hyperplane through the $n$ given points.

In terms of secant varieties, the varieties $\sigma_k(V^{2,n})$ can be identified
with the varieties of symmetric matrices of rank $\le k$ of order $(n+1)\times
(n+1)$, which have codimension
${{n-k+2}\choose 2}$.

The cases $d=4$, $2\le n\le 4$, $k={{n+2}\choose 2}-1$ are exceptional because
there is a (unique and smooth) quadric through the points, and the double
quadric is a quartic singular at the given points, while ${{n+4}\choose 4}\le
(n+1)\left[{{n+2}\choose 2}-1\right]$ exactly for $2\le n\le 4$.

The corresponding defective secant varieties $\sigma_k(V^{4,n})$ (with
$k={{n+2}\choose 2}-1$)
are hypersurfaces whose equation can be described as follows.

For any $\phi\in S^4V$, let $A_{\phi}\colon S^2V^{\vee}\to S^2V$ be the
contraction operator.
It is easy to check that if $\phi\in V^{4,n}$ then $rk A_{\phi}=1$ (by
identifying the Veronese variety with its affine cone).
It follows that if $\phi\in \sigma_k(V^{4,n})$ then $rk A_{\phi}\le k$. When
$k={{n+2}\choose 2}-1$ also the converse holds and $\det A_{\phi}=0$ is the
equation of the corresponding secant variety $\sigma_k(V^{4,n})$.
When $n=2$ the quartics in $\sigma_5(V^{4,2})$ are sum of five 4-powers of
linear forms and they are called Clebsch quartics \cite{Cl}.

The case $n=4,d=3,k=7$ is more subtle. In this case, since ${{7}\choose
3}=7\cdot 5$, it is expected that
no cubics exist with seven given singular points.
But indeed through seven points there is a rational normal curve $C_4$, which, in
a convenient system of coordinates, has equation
$$rk \left[\begin{array}{ccc}
x_0&x_1&x_2\\
x_1&x_2&x_3\\
x_2&x_3&x_4\\\end{array}\right]\le 1.$$
Its secant variety is the cubic with equation
$$\det \left[\begin{array}{ccc}
x_0&x_1&x_2\\
x_1&x_2&x_3\\
x_2&x_3&x_4\\\end{array}\right] = 0,$$
which is singular along the whole $C_4$. This is the same $J$ invariant
which describes harmonic $4$-ples on the projective line.
The paper \cite{CH} contains a readable  proof of the uniqueness
of the cubic singular along $C_4$.

Let us mention that in \cite{Rei} Reichstein gives an algorithm
to find if $f\in S^3({\bf C}^5)$ belongs to the hypersurface
$\sigma_7(V^{3,4})$. For the invariant
equation of
this hypersurface, which has degree $15$, see \cite{Ot}.

\section{Terracini's inductive argument}
\label{sezione4}

Terracini in \cite{Te3} considers a union $X$ of double points on ${\bf P}^3$
and studies the dimension of the system of hypersurfaces through $X$ by specializing some of the points to a
plane
${\bf P}^2\subseteq{\bf P}^3$. This is the core of an inductive procedure which has been considered
from several authors since then. The appealing fact of the inductive procedure is that
it covers almost all the cases with a very simple argument. This is the point that we want to explain in this section.
The remaining cases, which are left out because they do not fit the arithmetic of the problem, have to
be considered with a clever degeneration argument, which we postpone to Section~\ref{horace}.

Let $X$ be a  union of $k$ double points of ${\bf P}^{n}$, let ${\cal I}_X$ be the corresponding ideal sheaf
and fix a hyperplane
$H\subset {\bf P}^{n}$.
The trace of $X$ with respect to $H$ is the scheme $X\cap H$ and
the residual of $X$ is
the scheme $\widetilde{X}$ with ideal sheaf ${\cal I}_X:{\cal O}_{{\bf P}^n}(-H)$.
In particular if we specialize $u\leq
k$ points on the hyperplane $H$, the trace
$X\cap H$ is given by $u$ double points of ${\bf P}^{n-1}$, and
the residual $\widetilde{X}$ is given by $k-u$ double points and by $u$ simple points.

\begin{thm0}\label{caste}
Let $X$ be a union of $k$ double points of ${\bf P}^{n}$ and fix a hyperplane
$H\subset {\bf P}^{n}$ containing $u$ of them.
Assume that $X\cap H$ does impose independent conditions on ${\cal O}_H(d)$
 and  the residual $\widetilde{X}$ does impose independent
conditions on ${\cal O}_{{\bf P}^n}(d-1)$.
Assume moreover one of the following pair of inequalities:

(i) $un\le{{d+n-1}\choose{n-1}}\qquad  k(n+1)-un\le{{d+n-1}\choose{n}}$,

(ii) $un\ge{{d+n-1}\choose{n-1}}\qquad  k(n+1)-un\ge{{d+n-1}\choose{n}}$.

\noindent  Then  $X$ does impose independent conditions on the system ${\cal O}_{{\bf P}^n}(d)$.
\end{thm0}
{\it Proof.}
We want to prove that $I_X(d)$ has the expected dimension
$$\max\left({{d+n}\choose{n}}-k(n+1),0\right).$$
Taking the global sections of the restriction exact sequence
$$0\rig{} {\cal I}_{\widetilde{X}}(d-1)\rig{}{\cal I}_X(d) \rig{} {\cal I}_{X\cap H}(d)\rig{} 0,$$
we obtain the so called Castelnuovo exact sequence
\begin{equation}\label{castelnuovo}0\rig{}I_{\widetilde{X}}(d-1)\rig{}I_X(d)\rig{}I_{X\cap H}(d)
\end{equation}
from which we get the following inequality
$$\dim I_X(d)\leq \dim I_{\widetilde{X}}(d-1) + \dim I_{X\cap H}(d).$$
Since
$X\cap H$ imposes independent conditions on ${\cal O}_H(d)$
we know that $\dim I_{X\cap H}(d)=\max({{d+n-1}\choose{n-1}}-un,0)$;
on the other hand, since $\widetilde{X}$
imposes independent conditions on ${\cal O}_{{\bf P}^n}(d-1)$
it follows that
$\dim I_{\widetilde{X}}(d-1) =\max({{d-1+n}\choose{n}}-(k-u)(n+1)-u,0).$

Then in case {(i)}, we get $\dim I_X(d)\leq {{d+n}\choose{n}}-k(n+1)$,
while in case {(ii)}, we get $\dim I_X(d)\leq 0$.
But since $\dim I_X(d)$ is always greater or equal than the expected dimension, we conclude.
\qedd

In many cases a standard application of the above theorem gives  most of the cases of
Theorem \ref{AH}.

Let us see some examples in ${\bf P}^3$.
It is easy to check directly that there are no cubic surfaces with five singular points
(e.g. by choosing the five fundamental points in ${\bf P}^3$). This is the starting point of the
induction.

Now consider $d=4$ and a union $X$ of $8$ general double points. Setting $u=4$ we check that the
inequalities of case {(i)} of \thmref{caste} are satisfied.
Hence we specialize $4$ points on a hyperplane $H$ in such a way that they are general on $H$, then
by Theorem \ref{pidue} it follows that the trace $X\cap H$ imposes independent conditions on quartics.
On the other hand we consider the residual $\widetilde{X}$, given by $4$ double points outside $H$
and $4$ simple points on $H$. We know that the scheme $\widetilde{X}$ imposes independent conditions on cubics, since
the previous step implies that $4$ general double points do, and moreover we can add $4$ simple
points contained in a plane. This is possible because there exists no cubics
which are unions of a plane and a quadric through $4$ general double points.
\thmref{caste} applies and we conclude that $8$ general double points impose independent conditions
on ${\cal O}_{\P3}(4)$.
Notice that $9$ double points (one of the exceptional cases in \thmref{AH}) do not impose independent conditions on quartic surfaces. Indeed if we apply the same
argument we get as trace $5$ double points on ${\bf P}^2$, which do not impose independent conditions on quartics by
Theorem \ref{pidue}.

Consider now the case $d=5$. To prove that a general union of $14$
double points in ${\bf P}^3$ imposes independent conditions on quintics, it is enough
to specialize $u=7$ points on a plane in such a way that the trace is general and we apply induction.
On the other hand, also the residual imposes independent conditions
on quartics by induction and since there is no quartics
which are unions of a plane and a cubic through $7$ general double points.
Again \thmref{caste} applies and we can conclude  that
any collection of general double points imposes independent conditions
on ${\cal O}_{\P3}(5)$.

For $d\geq6$ we can apply this simple argument and by induction it is possible to prove that
 $k$ double points impose independent conditions on surfaces of degree $d$
with the following possible exceptions, for  $6\le d\le 30$:
$$(d,k)=(6,21), (9,55), (12,114), (15,204), (21, 506), (27, 1015), (30, 1364).$$

In particular if $d\neq 0\ {\mathrm{mod}}\ 3$, then it turns out that $k$ double points impose independent conditions
on surfaces of degree $d$.
To extend the result to the case   $d= 0\ {\mathrm{mod}}\ 3$ and the only possibly missing values of $k$
(that is $k=\lceil \frac{(d+3)(d+2)(d+1)}{24}\rceil$) is
much more difficult. We will do this job in full generality in Section $6$.

\section{The case of cubics}
\label{sezione5}

The inductive procedure of the previous section does not work
with cubics ($d=3$) because by restricting to a hyperplane we reduce to
quadrics which have defective behavior.
Nevertheless the case of cubics is the starting point
of the induction, so it is crucial. Alexander and Hirschowitz solved
this case in \cite{AH1}, by a subtle blowing up and by applying
the differential Horace's method (see Section 6).
Chandler solved this case with more elementary techniques in \cite{Ch2}.
In this section we give a shorter (and still elementary) proof.

Given $n$, we denote $k_n=\lfloor\frac{(n+3)(n+2)}{6}\rfloor$ and
$\delta_n={{n+3}\choose 3}-(n+1)k_n$. Notice that
$k_n=\frac{(n+3)(n+2)}{6}$ for $n\neq 2\ \mathrm{mod}\ 3$. If
$n=3p+2$,  we get $k_n=\frac{(n+3)(n+2)}{6}-\frac{1}{3}=\frac{(n+4)(n+1)}{6}$ and
$\delta_n=p+1=\frac{n+1}{3}$.

This simple arithmetic remark shows that the restriction to codimension three linear subspaces has the
advantage to avoid the arithmetic problems, and this is our new main idea.
In this section we will prove the following theorem, which immediately implies the case $d=3$ of \thmref{AH}.

\begin{thm0}\label{cubiche}
Let $n\neq 2\ \mathrm{mod}\ 3$, $n\neq 4$.
Then $k_n$ double points impose independent conditions on cubics.

Let $n=3p+2$, then $k_n$ double points and a zero dimensional scheme of length
$\delta_n$ impose independent conditions on cubics.
\end{thm0}

The proof of Theorem \ref{cubiche} relies on the following description.
\begin{prop0}\label{proprep}
Let $n\ge 5$ and let $L, M, N \subset\P n$ be general subspaces of
codimension $3$. Let $l_i$ (resp. $m_i, n_i$) with $i=1, 2, 3$ be
three general points on $ L$, (resp. $ M, N$). Then there are no
cubic hypersurfaces in $\P n$ which contain $ L\cup M\cup{ N}$ and
which are singular at the nine points $l_i, m_i, n_i$, with $i=1, 2,
3$.
\end{prop0}

{\it Proof.}
For $n=5,6,7$ it is an explicit computation, which can be easily performed
with the help of a computer. Indeed in
$\P 5$ it is easy to check that $I_{L\cup N\cup M, \P 5}(3)$ has
dimension $26$. Choosing three general points on each
subspace and imposing them as singular points for the cubics, one can
check that they impose $26$ independent conditions.
Analogously in the cases of $\P 6$ and $\P 7$
it is possible to compute $\dim I_{L\cup N\cup M, \P
  6}(3)=27$  and $\dim I_{L\cup N\cup M, \P 7}(3)=27$, and
imposing the nine singular points, one can check that we get
$27$ independent conditions.

For $n\ge 8$ the statement follows by induction on $n$.
Indeed if $n\ge 8$ it is easy to check that there are no quadrics
containing ${L}\cup{M}\cup{N}$.
Then given a general hyperplane
$H\subset \P n$ the Castelnuovo sequence induces the isomorphism
$$0\rig{}I_{{L}\cup{ M}\cup{ N},{\P n}}(3)\rig{}I_{\left({ L}\cup{ M}\cup{ N}\right)\cap H,H}(3)\rig{} 0$$
hence specializing the nine points on the hyperplane $H$, since the space
$I_{{ L}\cup{ M}\cup{ N},{\P n}}(2)$ is empty, we get
$$0 \rig{}I_{X\cup { L}\cup{ M}\cup{ N},{\P n}}(3)\rig{}I_{\left(X\cup
    { L}\cup{ M}\cup{ N}\right)\cap H,H}(3)$$
where $X$ denotes the union of the nine double points
supported at $l_i, m_i, n_i$ with $i=1, 2, 3$.
Then our statement immediately follows by induction.
\qedd

{\bf Remark.} Notice that \propref{proprep} is false for $n=4$. Indeed
$I_{L\cup N\cup M, \P 4}(3)$ has
dimension $23$ and there is a unique cubic singular at  the nine points
$l_i$, $m_i$, $n_i$, $i=1,2,3$.  Also the following \propref{proprep2} and \propref{proprep3}
are false
for $n= 4$, indeed their statements reduce to the statement of
Theorem \ref{cubiche}, because a cubic singular at $p$ and $q$ must contain
the line $<p,q>$ .

\begin{prop0}\label{proprep2}
Let $n\ge 3$, $n\neq 4$ and let ${L, M}\subset\P n$ be subspaces of codimension three.
Let $l_i$ (resp. $m_i$) with $i=1, \ldots n-2$ be  general points on ${L}$
(resp. ${M}$). Then there are no cubic hypersurfaces in $\P n$
containing ${L}\cup{M}$ which are singular
at the $2n-4$ points $l_i, m_i$ with $i=1,  \ldots n-2$ and at three general points
$p_i\in \P n$, with $i=1,2,3$.
\end{prop0}

{\it Proof.}
The case $n=3$ is easy and it was checked in Section 4.
For $n=5,7$ it is an explicit computation.
Indeed it is easy to check that $\dim I_{L\cup M, \P 5}(3)=36$
and that the
union of three general points on $L$, three general points on $M$ and three general points on $\P 5$
imposes $36$ independent conditions on the system $I_{L\cup M, \P 5}(3)$.
In the case $n=7$ one can easily check that
$\dim I_{L\cup M, \P 7}(3)=54$,  and that
the union of five general points on $L$, five general points on $M$ and three general points on $\P 5$
imposes $54$ independent conditions.

For $n=6$ or $n\ge 8$, the statement follows by induction from $n-3$ to $n$.
Indeed given a third general codimension three subspace ${N}$, we get the exact sequence
$$0\rig{}I_{{ L}\cup{ M}\cup{ N},{\P n}}(3)\rig{}
I_{{ L}\cup{M},{\P n}}(3)\rig{}I_{\left({ L}\cup{ M}\right)\cap { N},{ N}}(3)\rig{}0$$
where the dimensions of the three spaces in the sequence are respectively $27$,  $9(n-1)$ and $9(n-4)$.

Let $X$ denote the union of the double points supported at $p_1,p_2,p_3$, $l_i$ and $m_i$ with $i=1,\ldots,n-2$.
Let us specialize $n-5$ of the points $l_i$ (lying on
${L}$) to $L\cap N$, $n-5$ of the points $m_i$ (lying on ${M}$) to $M\cap N$ and the three points
$p_1,p_2,p_3$ to $N$. Then we obtain a sequence
$$0\rig{}I_{X\cup {L}\cup{ M}\cup{ N},{\P n}}(3)\rig{}
I_{X\cup { L}\cup{M},{\P n}}(3)\rig{}I_{\left(X\cup{ L}\cup{ M}\right)\cap { N},{ N}}(3)$$
where the trace $\left(X\cup{ L}\cup{ M}\right)\cap { N}$
 satisfies the assumptions on $N=\P{{n-3}}$ and we can apply induction.
Then we conclude, since the residual satisfies the hypotheses of
\propref{proprep}.\qedd

\begin{prop0}\label{proprep3}
Let $n\ge 3$, $n\neq 4$ and let ${L}\subset\P n$ be  a subspace of codimension three.

(i) If $n\neq 2\ \mathrm{mod}\ 3$ then there are no cubic hypersurfaces in $\P n$
which contain ${L}$ and which are singular at $\frac{n(n-1)}{6}$
general points $l_i$ on ${L}$ and at $(n+1)$ general points
$p_i\in \P n$.

(ii) If  $n= 2\ \mathrm{mod}\ 3$ then there are no cubic hypersurfaces in $\P n$
which contain ${L}$, which are singular at $\frac{(n+1)(n-2)}{6}$
general points $l_i$ on ${L}$  and at $(n+1)$ general points $p_i\in\P
n$, and which contain a general scheme $\eta$ supported at $q\in L$
such that ${\rm length}(\eta)=\delta_n$ and
${\rm length}(\eta\cap L)=\delta_n-1$.
\end{prop0}

{\it Proof.}
The case $n=3$ is easy and already checked in Section 4.
For $n=5$ let $e_i$ for $i=0,\ldots, 5$ be a basis of $V$ and choose
$L$ spanned by $p_i=[e_i]$ for $i=0, 1, 2$.
Consider the system of cubics  with singular points at $p_i$ for
$i=0,\ldots ,5$,
at $[e_0+\ldots +e_5]$ and at other two random points.
Moreover impose that the cubics of the system contain
a general scheme of length 2 supported at $[e_0+e_1+e_2]$. Note that
such cubics contain $L$.
A direct computation shows that this system is empty, as we wanted.
For $n=7$ the statement (i) can be checked, with the help of a computer, by
computing the tangent spaces to $V^{3,7}$ at seven general points of $L$ and
at eight general points. The condition that the cubic contains $L$ can be imposed by
another simple point on $L$.

For $n=6$ or $n\ge 8$ the statement follows by induction, and by the sequence
$$0\rig{}I_{{L}\cup{M},{\P n}}(3)\rig{}
I_{{L},{\P n}}(3)\rig{}I_{{L}\cap {M},{M}}(3)\rig{}0$$
where $M$ is a general codimension three subspace. Denoting by $X$
the union of the double points supported at
the points $l_i$ and $p_i$ (and of the scheme $\eta$ in case (ii)), we get
$$0\rig{}I_{X\cup{L}\cup{M},{\P n}}(3)\rig{}
I_{X\cup {L},{\P n}}(3)\rig{}I_{(X\cup{L})\cap {M},{M}}(3)$$

Assume now $n\neq 2\ \mathrm{mod}\ 3$. We specialize $\frac{(n-3)(n-4)}{6}$ of
the points $l_i$ to $L\cap M$ and $n-2$ of the points $p_i$ to $M$.
Thus we have left $n-2$ points general on $L$ and $3$ points general
on $\P n$ and we can use \propref{proprep2} on the residual and the induction on the trace.

If $n= 2\ \mathrm{mod}\ 3$,  we specialize $\frac{(n-2)(n-5)}{6}$ of the points
$l_i$ to $M\cap L$, and $n-2$ of the points $p_i$ and the scheme $\eta$
 to $M$ in such a way that $\eta\subset M$ 
and ${\rm length}(\eta\cap L)={\rm length}(\eta\cap L\cap M)=\delta_n-1$
(we can do this since $n\ge 8$)
and we conclude analogously.
\qedd

{\it Proof of Theorem \ref{cubiche}.}
We fix  a codimension three linear subspace $L\subset{\bf P}^n$ and we
prove the statement by induction by using the exact sequence
$$0\rig{}I_{{L},{\P n}}(3)\rig{}I_{\P n}(3)\rig{}I_L(3)$$

Assume first $n\neq 2\ \mathrm{mod}\ 3$.
We specialize to $L$ as many points as possible in order that
the trace with respect to $L$ imposes independent conditions on the cubics of
$L$. Precisely, we have $k_n=\frac{(n+3)(n+2)}{6}$ double
points and we  specialize $\frac{n(n-1)}{6}$ of them on $L$, leaving
$(n+1)$ points outside. Then the result follows from \propref{proprep3} and
by induction on $n$.
The starting points of the induction are $n=3$ (see Section $4$) and
$n=7$
(in this case it is enough to check that $15$ general tangent spaces to
$V^{3,7}$ are independent; notice that for $n=4$ the statement is false, see Section $3$).

In the case $n=2\ \mathrm{mod}\ 3$, we specialize
$k_{n-3}=\lfloor\frac{n(n-1)}{6}\rfloor=\frac{(n+1)(n-2)}{6}$
double points on $L$ and we leave $k_n-k_{n-3}=n+1$ double points outside $L$.
Moreover we specialize the scheme $\eta$
on $L$ in such a way that $\eta\cap L$ has length $\delta_n-1=\delta_{n-3}$.
Thus \propref{proprep3} applies again and we conclude by induction.
The starting point of the induction is $n=2$ (see \thmref{pidue}).\qedd

\section{The \ degeneration \ argument: \  ``la \ m\'ethode d'Horace dif\-f\'erentielle''}
\label{horace}
This section is devoted to the proof of Theorem \ref{AH} in the case $d\geq4$.

In order to solve the arithmetic problems revealed in the Section 4, Alexander and Hirschowitz
have introduced a clever degeneration argument, called the differential Horace's method.
We follow in this section the simplified version of the method performed by Chandler in
\cite{Ch}, trying to supply more details.  For the convenience of the reader
we describe first the case of sextics in $\P3$ (see \propref{ah36}), which is enough to understand the main idea.
In fact the pair $(6,21)$ was the first gap we met at the end of Section $4$.
After this case we will provide the proof in full generality.

Let $X, Z\subseteq{\bf P}^n={\bf P}(V)$ be  zero
dimensional subschemes, ${\cal I}_X$ and ${\cal I}_Z$ the corresponding ideal sheaves
and ${\cal D}={\cal I}_Z(d)$ for some $d\in {\bf N}$. The space ${\rm H}^0({\cal D})$
defines a linear system. The Hilbert function of $X$ with respect to ${\cal D}$ is defined as follows:
$$h_{{\bf P}^n}(X,{\cal D}):=\dim{\rm H}^0({\cal D})-\dim{\rm H}^0({\cal I}_X\otimes{\cal D}).$$
Notice that if ${\cal D}={\cal O}_{{\bf P}^n}(d)$,
then ${\rm H}^0({\cal I}_X\otimes{\cal D})=I_X(d)\subseteq
S^d V^{\vee}$  and we get
$$h_{{\bf P}^n}(X,d):=h_{{\bf P}^n}(X,{\cal O}(d))={{d+n}\choose{n}}-\dim I_X(d).$$
In other words $h_{{\bf P}^n}(X,d)$ is the codimension of the subspace $I_X(d)$ in the space
of homogeneous polynomials of degree $d$.

We say that $X$ {\em imposes independent conditions on ${\cal D}$} if
$$h_{{\bf P}^n}(X,{\cal D})=\min \left(\deg X, h^0({\cal D})\right)$$
This generalizes the definition given in Section 2 where ${\cal D}={\cal O}(d)$.

In particular if $h_{{\bf P}^n}(X,{\cal D})=\deg X$, we say that $X$ is  ${\cal
D}$-independent, and in the case  ${\cal D}={\cal O}(d)$, we say
$d$-independent.
Notice that if $Y\subseteq X$, then
if $X$ is ${\cal D}$-independent, then so is $Y$.
On the other hand if $h_{{\bf P}^n}(Y,d)={{d+n}\choose{n}}$, then
$h_{{\bf P}^n}(X,d)={{d+n}\choose{n}}$.

A zero dimensional scheme is called {\em curvilinear} if it is contained
in a non singular curve. A curvilinear scheme contained in a union
of $k$ double points has degree smaller than or equal to $2k$.

The following crucial lemma is due to Chandler \cite[Lemma 4]{Ch}.
\begin{lemma0}[Curvilinear Lemma]
\label{curvilinear}
Let $X\subseteq{\bf P}^n$ be a zero dimensional scheme
contained in a finite union of double points and $\cal D$ a linear system
on ${\bf P}^n$. Then $X$ is $\cal D$-independent if and only if
every curvilinear subscheme of $X$ is $\cal D$-independent.
\end{lemma0}
{\it Proof.}\ One implication is trivial. So let assume that
every curvilinear subscheme of $X$ is $\cal D$-independent.
Suppose first that $X$ is supported at one point $p$. We
prove the statement by induction on $\deg X$. If $\deg X=2$, then
$X$ is curvilinear and the claim holds true.

Now suppose $\deg X>2$ and let us prove that $h(X,{\cal D})=\deg
X$. Consider a subscheme $Y\subset X$ with $\deg Y=\deg X-1$. We
have $$h(Y,{\cal D})\leq h(X,{\cal D})\leq h(Y,{\cal D})+1.$$ By
induction $h(Y,{\cal D})=\deg Y=\deg X-1$. Then it is sufficient
to construct a subscheme $Y\subset X$ with $\deg Y=\deg X-1$
and $h(X,{\cal D})=h( Y,{\cal D})+1$.

In order to do this, consider a curvilinear subscheme $\xi\subset
X$, i.e.\ a degree $2$ subscheme of a double point. By hypothesis
we know that $\xi$ is ${\cal D}$-independent, i.e.\ $h(\xi,{\cal
D})=2$. Obviously we also have $h(p,{\cal D})=1$, where $p$
denotes the simple point. It follows that there exists a section
$s$ of ${\cal D}$ vanishing on $p$, and not on  $\xi$.
We define then $Y=X\cap Z$, where $Z$ is the zero locus of $s$.
Since $X$ is contained in a union of double points,
by imposing the condition $s=0$ we obtain
$\deg Y=\deg X-1$. Moreover $h(X,{\cal D})> h(Y,{\cal D})$ because
$s$ vanishes on $Y$ and does not on $X$. Then we conclude that
$$h(X,{\cal D})= h(Y,{\cal D})+1=\deg Y+1=\deg X.$$

Now consider $X$ supported at $p_1,\ldots,p_k$. Suppose by
induction on $k$ that the claim holds true for schemes supported
at $k-1$ points and we prove that $h(X,{\cal D})=\deg X$.
Let
$$A=X\cap{p_k^2}\quad\textrm{ and }\quad B=X\cap\{p_1,\ldots,p_{k-1}\}^2,$$
where $\{p_1,\ldots,p_{k-1}\}^2$ denotes the union of the double points $p_i^2$
and $X$ is a disjoint union of $A$ and $B$. Consider  ${\cal D}'={\cal I}_{B}\otimes{\cal D}$.

Let $\zeta$ be any curvilinear subscheme of $A$ and ${\cal
D}''={\cal D}\otimes {\cal I}_\zeta$. For every curvilinear
$\eta\subset B$ we have
$$h(\eta,{\cal D}'')=\dim {\rm H}^0({\cal D}\otimes {\cal I}_\zeta)
-\dim{\rm H}^0(I_{\zeta\cup\eta}\otimes{\cal D})=$$
$$=\dim{\rm H}^0({\cal D})-\dim{\rm H}^0({\cal I}_{\zeta\cup\eta}\otimes{\cal D})-\dim{\rm H}^0({\cal D})
+\dim{\rm H}^0({\cal D}\otimes {\cal I}_\zeta)=$$
$$=h(\zeta\cup\eta,{\cal D})-h(\zeta,{\cal D})=(\deg\zeta+\deg\eta)-\deg\zeta=\deg\eta$$
i.e.\ every curvilinear subscheme of $B$ is ${\cal
D}''$-independent. By induction it follows that $B$ is ${\cal
D}''$-independent, i.e.\ $h(B,{\cal D}\otimes {\cal I}_\zeta)=\deg B$.

Then we get  in the same way
$$h(\zeta\cup B,{\cal D})=h(\zeta,{\cal D})+h(B,{\cal D}'')=\deg\zeta+\deg B,$$
and again
$$h(\zeta,{\cal D}')=h(\zeta\cup B,{\cal D})-h(B,{\cal D})$$
hence putting together the last two equations and using the inductive assumption
we get
$$h(\zeta,{\cal D}') = (\deg\zeta+\deg B)-\deg B=\deg\zeta$$

We proved that every curvilinear subscheme of $A$ is ${\cal
D}'$-independent. Since $A$ is supported at one single point, from
the first part it follows that $A$ is ${\cal D}'$-independent.

Obviously ${\cal I}_{A}\otimes{\cal D}'={\cal I}_A\otimes {\cal
I}_B\otimes{\cal D}= {\cal I}_{X}\otimes{\cal D}.$ Then we
conclude, by using induction on $B$, that
$$h(X,{\cal D})=\dim{\rm H}^0({\cal D})-\dim{\rm H}^0({\cal I}_X\otimes{\cal D})=
\dim{\rm H}^0({\cal D})-\dim{\rm H}^0({\cal I}_A\otimes{\cal
D}')=$$
$$=h(B,{\cal D})+h(A,{\cal D}')=\deg B+\deg A=\deg X.$$
\qedd

Let us denote by $AH_{n,d}(k)$ the following statement:
{\it there exists a collection of $k$ double points
in ${\bf P}^n$ which impose
independent conditions on ${\cal O}_{{\bf P}^n}(d)$.}

Before considering the general inductive argument, we analyze in details
the first interesting example. We ask how many conditions
$21$ double points impose on ${\cal O}_{{\bf P}^3}(6)$ and we will prove that
$AH_{3,6}(21)$ holds true.

\begin{prop0} \label{ah36}
A collection of $21$ general double points imposes
independent conditions on ${\cal O}_{{\bf P}^3}(6)$.
\end{prop0}

{\it Proof.}
Notice that we cannot specialize $u$ points in such a way that conditions either (i)
or (ii) of Theorem \ref{caste}
are satisfied. Then we choose $u$ maximal such that
$nu<k(n+1)-{{d+n-1}\choose{n}}$, that is $u=9$.

By \thmref{pidue} and by Section \ref{sezione4} we know the following facts:

{\it (i)} $AH_{2,6}(9)$, and in particular $9$ general double points in $\bf{P}^2$ are $6$-independent,

{\it(ii)} $AH_{3,5}(12)$, and $12$ general double points in $\bf{P}^3$ are $5$-independent,

{\it(iii)} $AH_{3,4}(11)$, and there exist no quartic surfaces through $11$ general double points.

{\bf Step 1:}
Fix a plane ${\bf P}^{2}\subseteq{\bf
P}^{3}$. Let $\gamma\in{\bf P}^2$ be a point and $\Sigma$ a collection of $11$ general points
not contained in ${\bf P}^2$.
By {\it(ii)}, it follows that
$$h_{{\bf P}^{3}}(\{\gamma\}^2_{|{\bf P}^2}\cup\Sigma^2,5)=\deg (\{\gamma\}^2_{|{\bf P}^2}\cup\Sigma^2)=47.$$

{\bf Step 2:}
Now we want to add
a collection of $9$ points on ${\bf P}^2$
to the scheme $\{\gamma\}^2_{|{\bf P}^2}\cup\Sigma^2$.
It is obvious that if we add $9$ general simple points of ${\bf P}^3$
the resulting scheme would be $5$-independent.
But we want to add $9$ points contained in the plane. In fact we obtain the same conclusion
once we prove that there exists no quintic surface which is union of a plane and a quartic through
$\Sigma^2$. Indeed by  {\it{(iii)}} we know that
$\dim I_{\Sigma^2}(4)=0,$
hence we can choose a
collection $\Phi$ of $9$ simple points in ${\bf P}^2$ in such a way that
the scheme
${\{\gamma\}^2_{|{\bf P}^2}\cup\Sigma^2\cup\Phi}$ is $5$-independent.

{\bf Step 3:}
By {\it (i)}, it follows that
the scheme
$(\Phi^2_{|{\bf P}^{2}}\cup\gamma)\subseteq{\bf P}^2$ has Hilbert function
$$h_{{\bf P}^{2}}(\Phi^2_{|{\bf P}^{2}}\cup\gamma,6)=28$$ i.e.\ it
is $6$-independent.

\medskip

Now for $t\in{\bf K}$, let us choose a flat family
of general points
$\delta_t\subseteq{\bf P}^3$ and a
family of planes $\{H_t\}$ such that

$\bullet\quad \delta_t\in H_t$  for any $t$,

$\bullet\quad\delta_t\not\in{\bf P}^2$ for any $t\neq0$,

$\bullet\quad H_0={\bf
P}^2$ and $\delta_0=\gamma\in{\bf P}^2$.

\noindent Now consider the following schemes:
$\{\delta_t\}^2$,
${\Phi}^2$, where $\Phi$ is the collection of $9$ points introduced in Step $2$ and
$ {\Sigma}^2$, the collection of $11$ double points introduced in Step $1$.
Then in order to prove that $AH_{3,6}(21)$ holds,
 it is enough to prove the following claim.

{\bf Claim:} There exists $t\neq0$ such that the scheme $\{\delta_t\}^2$
is independent with respect to the system
$I_{\Phi^2\cup\Sigma^2}(6)$.

{\it Proof of the claim.}
Assume by contradiction that the claim is false. Then by Lemma
\ref{curvilinear} for all $t$ there exist pairs
$(\delta_t,\eta_t)$ with $\eta_t$
a curvilinear scheme supported in $\delta_t$ and contained in
$\{\delta_t\}^2$ such that
$$h_{{\bf P}^3}(\Phi^2\cup\Sigma^2\cup\eta_t,6)<82.
$$
Let $\eta_0$ be the limit of $\eta_t$.

By the semicontinuity of the Hilbert function and by
the previous inequality we get
\begin{equation}\label{ipotesiassurda}
h_{{\bf P}^3}(\Phi^2\cup\Sigma^2\cup \eta_0,6)\leq h_{{\bf P}^3}
(\Phi^2\cup\Sigma^2\cup \eta_t,6)<82.
\end{equation}

Consider the following two possibilities
\begin{itemize}
\item[1)]$\eta_0\not\subset\bf{P}^2$.

By applying the Castelnuovo exact sequence to $\Sigma^2\cup\Phi^2\cup \eta_0$, and by
using Step 2 and Step 3, we obtain
$$h_{\bf{P}^3}(\Sigma^2\cup\Phi^2\cup \eta_0,6)\geq h_{\bf{P}^3}(\Sigma^2\cup\Phi\cup\widetilde \eta_0,5)+
h_{\bf{P}^2}((\Phi^2_{|\bf{P}^2}\cup \gamma),6)=$$
$$=54+28=82,$$
a contradiction with (\ref{ipotesiassurda}).

\item[2)]$\eta_0\subset\bf{P}^2$.
By the semicontinuity of the Hilbert function
there exists an open neighborhood $O$ of 0 such that for any $t\in
O$
$$h_{{\bf P}^3}(\Phi\cup\Sigma^2\cup \{\delta_t\}^2_{|H_t},5)\geq
h_{{\bf P}^3}(\Phi\cup\Sigma^2\cup \{\gamma\}^2_{|\bf{P}^2},5)=9+44+3=56$$
and the equality holds.
In particular the subscheme\\ $\Phi\cup\Sigma^2\cup \eta_0\subset
\Phi\cup\Sigma^2\cup \{\gamma\}^2_{|\bf{P}^2}$ is $5$-independent, then
$h_{{\bf P}^3}(\Phi\cup\Sigma^2\cup \eta_t,5)=9+44+2=55$ for all $t\in O$.

Hence for any $t\in O$, by applying again the Castelnuovo exact sequence, we get
$$h_{{\bf P}^3}(\Phi^2\cup\Sigma^2\cup \eta_t,6)\geq h_{{\bf P}^3}(\Phi\cup\Sigma^2\cup\eta_t,5)+
h_{{\bf P}^2}(\Phi^2_{|{\bf P}^2},6)=55+27=82$$
contradicting again the inequality (\ref{ipotesiassurda}) above.
\end{itemize}
This completes the proof of the proposition.
\qedd

{\bf Remark.} We want to comment ``why" the proof of \propref{ah36} works.
A double point in $\P 3$ has length $4$; specializing it on a plane
we get a trace of length $3$ and a residual of length $1$. Among the $21$ points,
$9$ points are specialized on the plane ${\bf P}^2$, and $11$ remain outside.
 After this process has been performed, the trace defines a subspace of codimension $27$
in ${\rm H}^0(\O_{{\bf P}^2}(6))\simeq{\bf K}^{28}$ and there is no more room in the trace to specialize
the last point on ${\bf P}^2$, nor there is room in the residual to keep it outside. Thanks to the degeneration argument,
called the differential Horace's method, the last
point $\{\gamma\}^2$ ``counts like" a point of length $1$ in the trace,
and there is room for it.
This single point in the trace, which allows to solve the problem, reminds us of the Roman legend of the Horaces.

\medskip

In Theorem \ref{main} below we describe the general inductive argument.
It could be not enough to specialize only one point $\gamma$, in general we need to specialize $\epsilon$
points, with $0\le\epsilon < n$ to be chosen.  We need the following
easy numerical lemma, proved by Chandler \cite{Ch} in a slightly different form.

\begin{lemma0}\label{numerics}
Fix the integers $2\le n, 4\le d, 0\le k\le \lceil\frac{1}{n+1}{{n+d}\choose{n}}\rceil$ and let $u\in{\bf Z}$, $0\leq \epsilon<n$
 such that
$nu+\epsilon=k(n+1)-{{n+d-1}\choose{n}}.$ Then we have

(i) $n\epsilon + u\leq{{n+d-2}\choose{n-1}}$;

(ii) ${{n+d-2}\choose{n}}\le (k-u-\epsilon)(n+1)$;

(iii) $k-u-\epsilon\ge n+1$, for $d=4$ and $n\ge 10$.
\end{lemma0}

{\it Proof.}
We have
$$u\le  \frac{1}{n}\left({{n+d}\choose{n}}+(n+1)- {{n+d-1}\choose{n}}\right)=
\frac{1}{n}{{n+d-1}\choose{n-1}}+\frac{n+1}{n}
$$
hence
$$ n\epsilon + u \le n(n-1)+\frac{1}{n}{{n+d-1}\choose{n-1}}+\frac{n+1}{n}$$
and the right hand side is smaller than or equal to ${{n+d-2}\choose{n-1}}$
except for \linebreak
$(n,d)=(3,4), (4,4), (5, 4)$.
In this cases the inequality (i)  can be checked directly.

The inequality (ii) follows from (i) and
from the definition of $u$ and $\epsilon$.

In order to prove (iii) let us remark that by definition of $u$ we get
$$\mbox{\small
$k-u-\epsilon= -\frac{k}{n}+\frac{1}{n}{{n+3}\choose{n}}-\frac{(n-1)\epsilon}{n}\ge
\frac{1}{n}\left(-\frac{1}{n+1}{{n+4}\choose{n}}-1+{{n+3}\choose{n}}-(n-1)^2\right)$}$$
and the right hand side is greater or equal than $n+1$ for $n\ge 10$.
\qedd

\begin{thm0}\label{main}
Fix the integers $2\le n, 4\le d, \lfloor\frac{1}{n+1}{{n+d}\choose{n}}\rfloor\le k\le
\lceil\frac{1}{n+1}{{n+d}\choose{n}}\rceil$ and let $u\in{\bf Z}$, $0\leq \epsilon<n$
 such that
$nu+\epsilon=k(n+1)-{{n+d-1}\choose{n}}.$
Assume that
$AH_{n-1,d}(u)$
$AH_{n,d-1}(k-u)$,
$AH_{n,d-2}(k-u-\epsilon)$,
 hold.
Then $AH_{n,d}(k)$ follows.
\end{thm0}
{\it Proof.}
We will construct a scheme $\Phi^2\cup\Sigma^2\cup \Delta_t^2$ of $k$ double points which imposes
independent conditions on ${\cal O}_{{\bf P}^n}(d)$.

{\bf Step 1:}
Choose a hyperplane ${\bf P}^{n-1}\subseteq{\bf
P}^{n}$. Let $\Gamma=\{\gamma^1,\ldots ,\gamma^{\epsilon}\}$ be a collection of $\epsilon$
general points contained in ${\bf P}^{n-1}$ and $\Sigma$ a collection of $k-u-\epsilon$ points
not contained in ${\bf P}^{n-1}$.
By induction we know that $AH_{n,d-1}(k-u)$ holds, then it follows
$$h_{{\bf P}^{n}}(\Gamma^2_{|{\bf P}^{n-1}}\cup\Sigma^2,d-1)=\min\left((n+1)(k-u)-\epsilon,
{{n+d-1}\choose{n}}\right).$$
From the definition of $\epsilon$ it follows that
${{n+d-1}\choose{n}}=(n+1)(k-u)-\epsilon+u$
and since $u\geq0$, we obtain
$$h_{{\bf P}^{n}}(\Gamma^2_{|{\bf P}^{n-1}}\cup\Sigma^2,d-1)=(n+1)(k-u)-\epsilon.$$

{\bf Step 2:}
Now we want to add
a collection of $u$ simple points in ${\bf P}^{n-1}$
 to the scheme $\Gamma^2_{|{\bf P}^{n-1}}\cup\Sigma^2$ and we want to obtain a $(d-1)$-independent
scheme.
Notice that from Step 1 it follows that
$\dim I_{\Gamma^2_{|{\bf P}^{n-1}}\cup\Sigma^2}(d-1)=u$.
Thus it is enough to prove that there exist no hypersurfaces of degree $d-1$ which are
unions of ${\bf P}^{n-1}$ and of a hypersurface of degree $d-2$ through $\Sigma^2$.
In fact by induction we know that $\dim I_{\Sigma^2}(d-2)=
\max(0,{{n+d-2}\choose{n}}-(k-u-\epsilon)(n+1))$ and this
dimension vanishes by (ii) of Lemma \ref{numerics}.

Then
it follows that we can choose a
collection $\Phi$ of $u$ simple points in ${\bf
P}^{n-1}$ in such a way that the scheme\\
${\Gamma^2_{|{\bf P}^{n-1}}\cup\Sigma^2\cup\Phi}$ is $(d-1)$-independent, i.e.\
$$h_{{\bf P}^{n}}(\Gamma^2_{|{\bf P}^{n-1}}\cup\Sigma^2\cup\Phi,d-1)=
(n+1)(k-u)-\epsilon+u={{n+d-1}\choose{n}}.$$

\medskip

Now we split the proof in two cases.

{\bf First case:} $k(n+1)\leq{{d+n}\choose{n}}$.

{\bf Step 3:}
The assumption $k(n+1)\leq{{d+n}\choose{n}}$ implies that
$k= \lfloor\frac{1}{n+1}{{n+d}\choose{n}}\rfloor$
and $nu+\epsilon\leq{{d+n-1}\choose{n-1}}$.

By induction we know that $AH_{n-1,d}(u)$ holds, hence
the scheme $(\Phi^2_{|{\bf P}^{n-1}}\cup\Gamma)\subseteq{\bf P}^{n-1}$ has Hilbert function
$$h_{{\bf P}^{n-1}}(\Phi^2_{|{\bf P}^{n-1}}\cup\Gamma,d)=
\min({nu+\epsilon,{{d+n-1}\choose{n-1}}})=nu+\epsilon,$$
that is the scheme is $d$-independent.

\medskip

Now for $(t_1,\ldots, t_\epsilon)\in{\bf K}^{\epsilon}$, let us choose a flat family
of general points
$\{\delta^1_{t_1},\dots,\delta_{t_\epsilon}^{\epsilon}\}\subseteq{\bf P}^n$ and a
family of hyperplanes $\{H_{t_1},\ldots, H_{t_\epsilon}\}$ such that

$\bullet\quad\delta_{t_i}^i\in H_{t_i}$ for any $i=1,\ldots,\epsilon$
and for any $t_i$,

$\bullet\quad\delta_{t_i}^i\not\in{\bf P}^{n-1}$ for any $t_i\neq0$ and
for any $i=1,\ldots,\epsilon$,

$\bullet\quad H_0={\bf
P}^{n-1}$ and $\delta_0^i=\gamma^i\in{\bf P}^{n-1}$ for any $i=1,\ldots,\epsilon$.

\noindent Now let us consider the following schemes:

$\bullet\quad\Delta_{(t_1,\ldots, t_\epsilon)}^2=\{\delta^1_{t_1},\dots,\delta_{t_\epsilon}^{\epsilon}\}^2$,
  notice that
$\Delta_{(0,\ldots, 0)}^2={\Gamma}^2$;

$\bullet\quad {\Phi}^2$, where $\Phi$ is the collection of $u$ points introduced in Step $2$;

$\bullet\quad {\Sigma}^2$, the collection of $k-u-\epsilon$ double points introduced in Step $1$.

In order to prove that there exists a collection of $k$ points in $\mathbf{P}^n$ which impose
independent conditions on $\mathcal{O}_{\mathbf{P}^n}(d)$, it is enough to prove
the following claim.

{\bf Claim:} There exists $(t_1,\ldots, t_\epsilon)$ such that the scheme $\Delta_{(t_1,\ldots, t_\epsilon)}^2$
is independent with respect to the system
$I_{\Phi^2\cup\Sigma^2}(d)$.

{\it Proof of the claim.}
Assume by contradiction that the claim is false. Then by Lemma
\ref{curvilinear} for all $(t_1,\ldots, t_\epsilon)$ there exist pairs
$(\delta_{t_i}^i,\eta_{t_i}^i)$ for $i=1,\ldots,\epsilon$, with $\eta_{t_i}^i$
a curvilinear scheme supported in $\delta_{t_i}^i$ and contained in
$\Delta_{(t_1,\ldots, t_\epsilon)}^2$ such that
\begin{equation}\label{assurdo}
h_{{\bf
P}^n}(\Phi^2\cup\Sigma^2\cup\eta_{t_1}^1\cup\ldots,\eta_{t_\epsilon}^{\epsilon},d)<
(n+1)(k-\epsilon)+2\epsilon.
\end{equation}
Let $\eta_0^i$ be the limit of $\eta_{t_i}^i$, for
$i=1,\ldots,\epsilon$.

Suppose that $\eta_0^i\not\subset{\bf P}^{n-1}$ for $i\in
F\subseteq\{1,\ldots,\epsilon\}$ and $\eta_0^i\subset{\bf
P}^{n-1}$ for $i\in G=\{1,\ldots,\epsilon\}\setminus F$.

Given $t\in {\bf K}$, let us denote $Z_t^F=\cup_{i\in F}(\eta^i_t)$ and
$Z_t^G=\cup_{i\in G}(\eta^i_t)$. Denote by $\widetilde{\eta_0^i}$
the residual of $\eta_0^i$ with respect to ${\bf P}^{n-1}$ and by
$f$ and $g$ the cardinalities respectively of $F$ and $G$.

By the semicontinuity of the Hilbert function and by
(\ref{assurdo}) we get
\begin{equation}\label{seconda}
h_{{\bf P}^n}(\Phi^2\cup\Sigma^2\cup Z^F_0\cup Z^G_t,d)\leq h_{{\bf P}^n}
(\Phi^2\cup\Sigma^2\cup Z^F_t\cup
Z^G_t,d)<(n+1)(k-\epsilon)+2\epsilon.
\end{equation}

On the other hand, by the semicontinuity of the Hilbert function
there exists an open neighborhood $O$ of 0 such that for any $t\in
O$
$$h_{{\bf P}^n}(\Phi\cup\Sigma^2\cup(\cup_{i\in F}\widetilde{\eta_0^i})\cup Z^G_t,d-1)
\geq h_{{\bf P}^n}(\Phi\cup\Sigma^2\cup(\cup_{i\in
F}\widetilde{\eta_0^i})\cup Z^G_0,d-1).$$
Since $\Phi\cup\Sigma^2\cup(\cup_{i\in
F}\widetilde{\eta_0^i})\cup Z^G_0
\subseteq \Phi\cup\Sigma^2\cup\Gamma^2_{|{\bf P}^{n-1}}$,
by Step 2 we compute
$$h_{{\bf P}^n}(\Phi\cup\Sigma^2\cup(\cup_{i\in
F}\widetilde{\eta_0^i})\cup Z^G_0,d-1)=u+(n+1)(k-u-\epsilon)+f+2g.$$
Since $\Phi^2_{|{\bf P}^{n-1}}\cup
(\cup_{i\in F}\gamma_i)$ is a subscheme of $\Phi^2_{|{\bf P}^{n-1}}\cup\Gamma$, by Step 3
it follows that
$$h_{{\bf P}^{n-1}}(\Phi^2_{|{\bf P}^{n-1}}\cup
(\cup_{i\in F}\gamma_i),d)\geq nu+f$$
Hence for any $t\in O$, by applying the Castelnuovo exact sequence
to the scheme $\widetilde\Phi\cup\Sigma\cup Z^F_0\cup
Z^G_t$, we get,
$$h_{{\bf P}^n}(\Phi^2\cup\Sigma^2\cup Z^F_0\cup Z^G_t,d)\geq$$
$$\geq h_{{\bf P}^n}(\Phi\cup\Sigma^2\cup(\cup_{i\in F}\widetilde{\eta_0^i})\cup Z^G_t,d-1)+
h_{{\bf P}^{n-1}}(\Phi^2_{|{\bf P}^{n-1}}\cup
(\cup_{i\in F}\gamma_i),d)\geq$$
$$\geq(u+(n+1)(k-u-\epsilon)+f+2g)+(nu+f)=$$
$$=(n+1)(k-\epsilon)+2\epsilon,$$
contradicting the inequality (\ref{seconda}) above.
This completes the proof of the claim and of the first case.

\medskip

{\bf Second case:} $k(n+1)>{{d+n}\choose{n}}$.

It follows that $k= \lceil\frac{1}{n+1}{{n+d}\choose{n}}\rceil$
and $nu+\epsilon>{{d+n-1}\choose{n-1}}$.

If  ${{d+n-1}\choose{n-1}}-nu<0$ then we are in the easy case (ii) of \thmref{caste}
(indeed the second inequality of (ii) is equivalent to $\epsilon\ge 0$). Then
$AH_{n,d}(k)$ holds by applying \thmref{caste}. Indeed the assumptions of \thmref{caste}
are satisfied: in particular the assumption on the trace
follows from $AH_{n-1,d}(u)$, while the assumption on the residual follows from
$AH_{n,d-1}(k-u)$,
and $AH_{n,d-2}(k-u-\epsilon)$, which in particular implies $AH_{n,d-2}(k-u)$
by Step 2.

So we may assume that $0\le \nu:={{d+n-1}\choose{n-1}}-nu<\epsilon$.

{\bf Step 3:}
Differently from the first case, now we obtain
$$h_{{\bf P}^{n-1}}(\Phi^2_{|{\bf P}^{n-1}}\cup\Gamma,d)={{d+n-1}\choose{n-1}}<nu+\epsilon.$$
Note that if we substitute to $\Gamma$ its subset
$\overline{\Gamma}=\{\gamma_1,\ldots ,\gamma_{\nu}\}$
we get
$$h_{{\bf P}^{n-1}}(\Phi^2_{|{\bf P}^{n-1}}\cup\overline{\Gamma},d)=
{{d+n-1}\choose{n-1}}=nu+\nu$$
and the advantage of this formulation is that now we can apply
Lemma \ref{curvilinear} to the scheme $\Phi^2_{|{\bf P}^{n-1}}\cup\overline{\Gamma}$.

\medskip

Now choose a flat family
of general points
$\{\delta^1_{t_1},\dots,\delta_{t_\epsilon}^{\epsilon}\}\subseteq{\bf P}^n$ and a
family of hyperplanes $\{H_{t_1},\ldots, H_{t_\epsilon}\}$ with the same properties as above.

Let us denote
$$\overline{\Delta}_{(t_1,\ldots,t_\epsilon)}=\{\delta^1_{t_1},\dots,\delta_{t_\nu}^{\nu}\}^2\cup
\{\delta_{t_{(\nu+1)}}^{\nu+1}\}^2_{|H_{t_{(\nu+1)}}}\cup\dots
\{\delta_{t_{\epsilon}}^{\epsilon}\}^2_{|H_{t_\epsilon}}.$$
Since obviously we have
$$h_{{\bf P}^n}(\Phi^2\cup\Sigma^2\cup \Delta_{(t_1,\ldots,t_\epsilon)}^2,d)\geq
h_{{\bf P}^n}(\Phi^2\cup\Sigma^2\cup \overline{\Delta}_{(t_1,\ldots,t_\epsilon)},d),$$
in order to conclude it is enough to prove the following claim.

{\bf Claim:}
There exists ${(t_1,\ldots,t_\epsilon)}$ such that the scheme $\overline{\Delta}_{(t_1,\ldots,t_\epsilon)}$
is independent with respect to the system
$I_{\Phi^2\cup\Sigma^2}(d)$.

We can prove the claim exactly as in the first case.
Indeed note that \linebreak
$\{\nu+1,\ldots ,\epsilon\}\subseteq G$.
Then $\Phi^2_{|{\bf P}^{n-1}}\cup
(\cup_{i\in F}\gamma_i)$ is a subscheme of
$\Phi^2_{|{\bf P}^{n-1}}\cup\overline{\Gamma}$, hence by \lemref{curvilinear}
it follows again that
$$h_{{\bf P}^{n-1}}(\Phi^2_{|{\bf P}^{n-1}}\cup
(\cup_{i\in F}\gamma_i),d)\geq nu+f$$

So the above proof of the claim works smoothly. This completes the proof of the second case.\qedd

\thmref{main} allows us to prove Theorem \ref{AH}, once we
have checked the initial steps of the induction.
Thanks to \thmref{cubiche}, the only problems occurring in the initial steps
depend on quadrics and on the exceptional cases. It is easy to see that the only cases we have to
study explicitly are
${\cal O}_{{\bf P}^n}(4)\textrm{ for }5\le n\le 9$.
Indeed for $n\ge 10$ we can apply (iii) of \lemref{numerics} and
the easy fact that $AH_{n,2}(k)$ holds if $k\ge n+1$,
because there are no quadrics with $n+1$ general double points.

Even for $n=9$ we have $k=71$ (or respectively $72$),  $(u,\epsilon)=(54,4)$,
(respectively $(55,5)$)
and still $k-u-\epsilon\ge n+1$  so that $AH_{n,2}(k-u-\epsilon)$ holds,
moreover we need $AH_{8,4}(54)$ (respectively $AH_{8,4}(55)$) and $AH_{9,3}(17)$
that will turn out to hold by the induction procedure.
The same argument applies for $n=6, 8$.

 For $n=7$, we have to consider $k=41$ or $42$.
 For $k=41$ it applies \thmref{caste} (i) with $u=30$, while for $k=42$ it applies
 \thmref{caste} (ii) again with $u=30$.

 In the remaining case $n=5$ we have $k=21$ and neither \thmref{caste} nor \thmref{main} apply
 because we always need $AH_{4,4}(14)$ which does not hold and indeed it is the last exceptional case of
 \thmref{AH}.
This case can be checked explicitly, by verifying that $21$ general tangent spaces
to $V^{4,5}$ are independent, with the help of a computer, or by an ad hoc argument, either as in
\cite{AH0} or as in the last paragraph of \cite{Ch}.

This completes the proof of Theorem \ref{AH}.
\medskip

{\bf Remark.} Alexander and Hirschowitz called the assumption $AH_{n-1,d}(u)$ in
\thmref{main}  the {\it dime} (lower dimension)   and the other assumptions
the {\it degue} (lower degree).

\section{Historical remarks}
\label{sezione7}
\subsection{The one dimensional case and the Sylvester Theorem}
In the case $n=1$ the Veronese variety $V^{d,1}$ is the rational normal curve $C_d$.
It is easy to check that
the higher secant variety $\sigma_k(C_d)$ has always the expected dimension
(moreover this is true for arbitrary curves, see \cite[Example V.1.6]{Z}).
In the setting of \thmref{AH} this follows from the fact that the space of one variable polynomials,
with given roots of fixed
multiplicities, has always the expected dimension.
Indeed there are well known explicit interpolation formulas
to handle this problem which go back to Newton and Lagrange.

The equations of the higher secant varieties to the rational normal curves $C_d$
were computed  by Sylvester in 1851.
In modern notation, given a vector space $U$ of dimension two and $\phi\in S^{2m}U$
it is defined the contraction operator $A_{\phi}\colon S^mU^{\vee}\rig{}S^mU$
and we have that $\phi\in \sigma_k(C_{2m})$ if and only if $rk A_{\phi}\le k$,
while in the odd case   we have $\phi\in S^{2m+1}U$,
the contraction operator $A_{\phi}\colon S^mU^{\vee}\rig{}S^{m+1}U$
and again we have that $\phi\in \sigma_k(C_{2m+1})$ if and only if $rk A_{\phi}\le k$.
It turns out that the equations of the higher secant varieties of the rational normal curve are given by
the minors of $A_{\phi}$. The matrices representing $A_{\phi}$ were called catalecticant by Sylvester \cite{Sy}.
In 1886 Gundelfinger (\cite{Gu})  treated the same problem from a different point of view by finding
the covariants defining $\sigma_k(C_{d})$ in the setting of classical invariant theory.
In \cite{Sy} Sylvester also found the canonical form of a general  $\phi\in S^{2m+1}U$
as sum of $m+1$ uniquely determined powers of linear forms. This is the first case of
the Waring problem for polynomials.

Making precise the  statement of Sylvester, we denote
$$f_{p,q}=\frac{\partial^{p+q}f}{\partial x^{p}\partial y^{q}}$$
and we get the following
\begin{thm0}[Sylvester]
\label{sylvester}
Let $f(x,y)$ be a binary form of degree $2m+1$ over the complex numbers.
Consider the $(m+1)\times (m+1)$ matrix $F$ whose $(i,j)$ entry
is $f_{2m-i-j,i+j}$ for
$0\le i,j\le m$ and denote $g(x,y)=\det F$.

(i) If $g(x,y)$ does vanish identically then $f\in\sigma_m(C_{2m+1})$,
and the converse holds.

(ii) If $g(x,y)$ does not vanish identically then
factorize
$$g(x,y)=\prod_{i=1}^{m+1}(p_ix+q_iy).$$
There are uniquely determined constants $c_i$ such that
$$f(x,y)=\sum_{i=1}^{m+1}c_i(p_ix+q_iy)^{2m+1}$$
if and only if  $g(x,y)$ has distinct roots.
(A convenient choice of $p_i, q_i$ allows of course to take $c_i=1$.)
\end{thm0}

It is worth to rewrite and prove Sylvester theorem
in the first nontrivial case, which is the case of quintics, as Sylvester himself did. The general case is analogous.
Let
$$f=a_0x^5+5a_1x^4y+10a_2x^3y^2+10a_3x^2y^3+5a_4xy^4+a_5y^5.$$
We have that $f\in \sigma_k(C_{5})$ if and only if
$$rk\left[\begin{array}{cccc}a_0&a_1&a_2&a_3\\
a_1&a_2&a_3&a_4\\
a_2&a_3&a_4&a_5\\
\end{array}\right]\le k.$$
We have  the formula
$$\frac{1}{5!}\left[\begin{array}{ccc}
f_{4,0}&f_{3,1}&f_{2,2}\\
f_{3,1}&f_{2,2}&f_{1,3}\\
f_{2,2}&f_{1,3}&f_{0,4}\\
 \end{array}\right] =
\left[\begin{array}{ccc}a_0x+a_1y&a_1x+a_2y&a_2x+a_3y\\
a_1x+a_2y&a_2x+a_3y&a_3x+a_4y\\
a_2x+a_3y&a_3x+a_4y&a_4x+a_5y\\
 \end{array}\right]$$
moreover Sylvester found the following equality between determinants
$$ \left|\begin{array}{ccc}a_0x+a_1y&a_1x+a_2y&a_2x+a_3y\\
a_1x+a_2y&a_2x+a_3y&a_3x+a_4y\\
a_2x+a_3y&a_3x+a_4y&a_4x+a_5y\\
 \end{array}\right|=
\left|\begin{array}{cccc}y^3&-x^2y&x^2y&-x^3\\
a_0&a_1&a_2&a_3\\
a_1&a_2&a_3&a_4\\
a_2&a_3&a_4&a_5\\
\end{array}\right|$$
and  Cayley pointed out to him (\cite{Sy}) that it follows from
$$\mbox{ \scriptsize $
\left[\begin{array}{cccc}y^3&-x^2y&x^2y&-x^3\\
a_0&a_1&a_2&a_3\\
a_1&a_2&a_3&a_4\\
a_2&a_3&a_4&a_5\\
\end{array}\right] \cdot
\left[\begin{array}{cccc}1&x&0&0\\
0&y&x&0\\
0&0&y&x\\
0&0&0&y\\
\end{array}\right] =
\left[\begin{array}{cccc}y^3&0&0&0\\
a_0&a_0x+a_1y&a_1x+a_2y&a_2x+a_3y\\
a_1&a_1x+a_2y&a_2x+a_3y&a_3x+a_4y\\
a_2&a_2x+a_3y&a_3x+a_4y&a_4x+a_5y\\
 \end{array}\right]. $} $$

We get that
$f\in \sigma_2(C_{5})$ if and only if
$$\left|\begin{array}{ccc}
f_{4,0}&f_{3,1}&f_{2,2}\\
f_{3,1}&f_{2,2}&f_{1,3}\\
f_{2,2}&f_{1,3}&f_{0,4}\\
 \end{array}\right| \equiv 0$$
(this is one of Gundelfinger's covariants)
and this proves (i).

In case (ii) we have the factorization
$$\left|\begin{array}{ccc}
f_{4,0}&f_{3,1}&f_{2,2}\\
f_{3,1}&f_{2,2}&f_{1,3}\\
f_{2,2}&f_{1,3}&f_{0,4}\\
 \end{array}\right| = (p_1x+q_1y)(p_2x+q_2y)(p_3x+q_3y)$$
and Sylvester proves in \cite{Sy} the ``remarkable discovery'' that there are  constants $c_i$ such that
$$f= c_1(p_1x+q_1y)^5+c_2(p_2x+q_2y)^5+c_3(p_3x+q_3y)^5$$
if and only if the three roots are distinct.

In particular the three linear forms $p_ix+q_iy$ are uniquely determined,
so that we get generically a canonical form as a sum of three 5-th powers.
The proof goes as follows.
Consider the covariant
$$g(a,x,y)= \left|\begin{array}{cccc}y^3&-x^2y&xy^2&-x^3\\
a_0&a_1&a_2&a_3\\
a_1&a_2&a_3&a_4\\
a_2&a_3&a_4&a_5\\
\end{array}\right|$$
which is called apolar to $f$ (we do not need this concept).
To any catalecticant matrix
$$  A = \left[\begin{array}{cccc}
a_0&a_1&a_2&a_3\\
a_1&a_2&a_3&a_4\\
a_2&a_3&a_4&a_5\\
\end{array}\right]    $$
such that $rk A=1$, it is associated a unique   $(x,y)\in{\bf P}^1$
such that
$$rk \left[\begin{array}{cccc}y^3&-x^2y&xy^2&-x^3\\
a_0&a_1&a_2&a_3\\
a_1&a_2&a_3&a_4\\
a_2&a_3&a_4&a_5\\
\end{array}\right]=1   $$
(it is easy to see this by looking at the parametric equations of the rational normal curve).

Assume now that the general catalecticant matrix $A$
is the sum of three catalecticant matrices of the same shape $A^i$ of rank $1$.
We may write $a=a^1+a^2+a^3$.
Let $(x_i,y_i) \in{\bf P}^1$ be the point associated to $a^i$.
Now compute
$g(a^1+a^2+a^3,x_1,y_1)$.
By linearity on rows, the determinant splits in $27$ summands,
among them there are $19$ which contain a row in $A^1$, which vanish because
any row of $A^1$ is dependent with   $(y_1^3,-x_1^2y_1,x_1y^2_1,-x_1^3)$,
and other $8$ which vanish because by the pigeon-hole principle they contain at least two rows from $A^2$ or
from $A^3$. It follows that  $g(a^1+a^2+a^3,x_1,y_1)=0$, then $(x_1,y_1)$ is a root of the covariant
$g(a,x,y)$. Since the same argument works also for $(x_i,y_i)$ with $i=2,3$, this ends the proof of
the uniqueness in Sylvester theorem.

To show the existence, we consider the $SL(U)$-equivariant morphism
$$\P{{}}(S^5U)\setminus \sigma_2(C_5) \rig{\pi} \P{{}}(S^3U)$$
defined by the covariant $g$.
The fiber of a polynomial $$z(x,y)=(p_1x+q_1y)(p_2x+q_2y)(p_3x+q_3y)\in\P{{}}(S^3U)$$ with distinct roots satisfies
\begin{equation}\label{fiber}
\mbox{\scriptsize
$\pi^{-1}(z)\supseteq\{c_1(p_1x+q_1y)^5+c_2(p_2x+q_2y)^5+c_3(p_3x+q_3y)^5| c_1\neq 0, c_2\neq 0, c_3\neq 0\}$
} \end{equation}
by the uniqueness argument and the fact that
 if some $c_i=0$ then the corresponding polynomial belongs to $\sigma_2(C_5)$.
Hence any  polynomial which is a sum of three distinct $5$-th powers must belong to one of the above fibers,
so that its image under $\pi$ must have three distinct roots.
Now a infinitesimal version of the above computation shows that if
$a=a^1+a^{11}+a^3$ where $a^{11}$ is on the tangent line at $a^1$,
then $g(a,x,y)$ has a double root at ($x_1,y_1)$.

In particular if $f\in\P{{}}(S^5U)$ cannot be expressed as the sum of three distinct $5$-th powers then
$\pi(f)$ must have a double root. This shows that the equality holds in (\ref{fiber})
and it concludes the proof.
\medskip

Note that the fiber of the general point is the algebraic torus given by
the $3$-secant $\P2$ minus three lines.
To make everything explicit,
denote by $T^i_p$ the $i$-th osculating space at $p$ to $C_5$,
so $T^1_p$ is the usual tangent line at $p$.
If $z(x,y)=(p_1x+q_1y)^2(p_2x+q_2y)\in\P{{}}(S^3U)$  then
$$\mbox{ \footnotesize
$\pi^{-1}(z)=<T^1_{(p_1x+q_1y)^5}, (p_2x+q_2y)^5>
\setminus \left( T^1_{(p_1x+q_1y)^5}\cup <(p_1x+q_1y)^5,(p_2x+q_2y)^5>\right)$ }$$
while if $z(x,y)=(p_1x+q_1y)^3\in\P{{}}(S^3U)$  then
$$\pi^{-1}(z)=T^2_{(p_1x+q_1y)^5}\setminus T^1_{(p_1x+q_1y)^5}$$
The last two fibers contain polynomials which can be expressed as sum of more than three powers.

In general  we consider the $SL(U)$-equivariant morphism
$$\P{{}}(S^{{2m+1}}U)\setminus \sigma_m(C_{2m+1}) \rig{\pi} \P{{}}(S^{m+1}U)$$
It follows that the polynomials $f\in\P{{}}(S^{{2m+1}}U)$ which have a unique canonical form as sum of $m+1$ powers
are exactly those lying outside the irreducible hypersurface which is the closure of $\pi^{-1}(\textrm{discriminant})$,
which has degree $2m(m+1)$, and it is the Zariski closure of the union of all
linear spans
$<T^1_{p_1}, p_2,\ldots , p_{m}>$
where $p_i$ are distinct points in $C_{2m+1}$.
If $z\in\P{{}}(S^{{m+1}}U)$ has $q$ distinct roots, then the fiber $\pi^{-1}(z)$ is isomorphic to $\P m$ minus
$q$ hyperplanes.

We emphasize that this argument by Sylvester not only proves the uniqueness of the canonical form of an odd
binary form as the sum of powers, but its also gives an algorithm to construct  it,
up to factor a polynomial equation in one variable.

A proof of \thmref{sylvester} using symbolic (umbral) calculus can be found in
\cite{KR}.

\subsection{The general case}
The cases of small degree and the first exceptions in Theorem \ref{AH}
were known since a long time.  The first nontrivial exception of plane quartics
was studied by Clebsch \cite{Cl}, who found in 1861 the equation of the degree $6$ invariant,
which gives the hypersurface $\sigma_5(\P 2, \O(4))$, as we sketched in Section \ref{sezione3}.  Richmond
in \cite{Ri} listed all the exceptions appearing in Theorem \ref{AH}.
For example in the more difficult case,
concerning a general cubic in $\P 4$ which is not the sum of seven cubes,
the method of Richmond is to construct the rational normal curve
through seven points, and then to manipulate
the equations of the problem into partial fractions.
A sentence from Richmond paper is illuminating:
{\it ``It does not appear to be possible to make any general application of the
method. I therefore continue to consider special problems''}.

To the best of our knowledge, the first paper which faces the problem (with $n\ge 2$) in general
was published by Campbell in 1892 \cite{Ca} on the ``Messenger of Mathematics'',  a journal which
stopped being published in 1928 and was absorbed by the Oxford Quarterly Journal.
Campbell is better known for the Campbell-Hausdorff formula for multiplication of exponents in Lie algebras.
He proved
an equivalent form of the second Terracini Lemma \ref{ter2} for linear systems of plane curves
by looking at the Jacobian of the system.
Campbell deduced that if a union $X$ of $k$ double points does not impose independent conditions
on plane curves of degree $d$, then  every  curve $C$ of degree $d$ through $X$ has to be a double curve, and $d$ is even.
The correct conclusion is that $C$ contains a double component, but it is easy to complete this argument,
as we saw in \thmref{pidue} and we repeat in a while.
The idea of Campbell was to add $t$ points in order that $3k+t={d+2\choose 2}$
and he found also the other equation $k+t={(d/2)+2\choose 2}-1$.
This system has only the two solutions
$$\left\{\begin{array}{c}d=2\\k=2\\t=0\end{array}\right.\quad\hbox{and}\quad
\left\{\begin{array}{c}d=4\\k=5\\t=0\end{array}\right.$$
which give the two exceptions of Theorem \ref{AH} for $n=2$.

Campbell then considered the case $n=3$ and he claimed  that if a
union $X$ of double points does not impose independent conditions
on surfaces of degree $d$, then  every  surface $C$ of
degree $d$ through $X$ has to be a double surface, and $d$ is
even. Although the conclusion is correct, the argument given by
Campbell seems to be wrong,
 otherwise it should work also when
$n=4$, but in this case the only cubic singular at seven points is actually reduced.
This fourth exceptional case
in the list of Theorem \ref{AH} was probably not known to Campbell.
It is worth to remark that Campbell proved
in the same paper that the only Veronese surfaces
which are weakly defective (in the modern notation, according to \cite{CC})
are given by the linear systems $|\O(d)|$ with $d=2, 4$ or $6$.  His argument is a slight modification
of the previous one, and it seems essentially correct.

Campbell concluded by applying his theorem to the canonical forms
of general hypersurfaces as sums of powers, and he got that the
expected number of summands is attained, with the only exceptions
of Theorem \ref{AH} (here $n\le 3$).
He did not apply Lasker Proposition \ref{lasker}. His more indirect
approach, which uses the Jacobian,
seems essentially equivalent to Proposition \ref{lasker}.

Campbell paper was not quoted by Richmond, we do not know if this is a signal of the
rivalry between Oxford and Cambridge.

In Italy the problem was faced in the same years by the school of
Corrado Segre. Palatini, a student of Segre, attacked the general
problem, and was probably not aware of Campbell's results. The
paper \cite{Pa1} is contemporary to \cite{Ri}, and treats the same
problem of the defectivity of the system of cubics in $\P 4$.
Palatini's argument that shows the defectivity is geometrical,
and resembles the one we have sketched in Section \ref{sezione3}.
A proof of Theorem \ref{AH} in the case $n=2$ is given in \cite{Pa2}.
We sketch the argument of Palatini in the case $d=7$,
which is direct, in opposition with the ones of Campbell and Terracini
which rely on infinitesimal computations. Palatini's aim is to prove that the
$12$-secant spaces to the $7$-Veronese embedding of $\P 2$ fill the ambient
space $\P {{35}}$. Denote by $D_p$ a plane curve of degree $p$.
Palatini first proved the following preliminary lemma.
\begin{lemma0}[Palatini]
\label{palatini}
(i) Assume  $p_1,\ldots ,p_{12}$ are general points in $\P 2$
and \linebreak
 $p_{13},\ldots ,p_{24}$ are chosen
such that $h^0(7H-\sum_{i=1}^{24} p_i)=36-24+1=13$
(one more than the expected value).
Then $p_1,\ldots ,p_{24}$ are the complete intersection of a $D_4$ with a $D_6$.

(ii) Conversely, if $Z=D_4\cap D_6$ then  $h^0(I_Z(7))=13$.
\end{lemma0}
{\it Proof.} By assumption a septic $D_7$ which contains $23$ of the given points
contains also the last one.
Let $D_3$ be the cubic through $p_1,\ldots ,p_{9}$.
Let $D_4$ be a quartic through $p_{11},\ldots ,p_{24}$;
by assumption it contains also $p_{10}$.
Considering the cubic through $p_1,\ldots ,p_8,p_{10}$,
it follows that $D_4$ contains also $p_9$, and continuing in this way, all the
points are contained in $D_4$.
The general sextic $D_6$ through $p_1,\ldots ,p_{24}$ does not contain $D_4$
as a component.
Indeed let $D_1$ be the line through $p_1$ and $p_2$. Let $D_6$ be a sextic
through $p_{4},\ldots ,p_{24}$,
by assumption it contains also $p_3$. Starting from other lines, such a $D_6$
contains all the $24$ points.
Then $H^0(6H-\sum_{i=1}^{24}p_i)= H^0(6H-\sum_{i=4}^{24}p_i)$ which has
dimension
$\ge 28-21=7>6=h^0(2H)$. This proves (i).
Part (ii) is today obvious from the Koszul complex. \qedd

By duality, a $12$-secant space $\pi$ corresponds
to the linear system of $D_7$ through  $12$ points $p_1,\ldots ,p_{12}$.
Consider all the other $12$-secant spaces which
meet our $\pi$. These correspond to collections of $12$ points
$p_{13},\ldots ,p_{24}$
such that $h^0(7H-\sum_{i=1}^{24} p_i)=13$. By Lemma  \ref{palatini}
these collections of $12$ points are parametrized by the pairs $(D_4,E)$ where
$D_4$ is a quartic through $p_1,\ldots ,p_{12}$ and $p_1+\ldots +p_{12}+E$ is a
divisor cut on $D_4$
 by a sextic.
There are $\infty^2$ quartic curves  and by Riemann-Roch formula
$E$ has $9$ parameters, so that there $\infty^{11}$  $12$-secant spaces which
meet our $\pi$. This means that for a general point of $\pi$ there are only
finitely many
$12$-secant spaces, hence the $12$-secant variety has the expected dimension as
we wanted.
Closing the paper \cite{Pa2}, Palatini wrote:
{\it ``si pu\`o gi\`a prevedere che l'impossibilit\`a di rappresentare
una forma s-aria generica con la somma di potenze di forme lineari
contenenti un numero di costanti non inferiore a quello contenuto
nella forma considerata, si avr\`a soltanto in casi particolari.''}
\footnote{\it One can expect that the impossibility of representing
a general form in $s$ variables as a sum of powers
of linear forms containing a number of constants not smaller
than the number of constants contained in the given forms,
holds only in a few particular cases. }
Then he listed the particular cases known to him,
and they are exactly the exceptions of Theorem \ref{AH}.
So this sentence can be considered as the first conjecture of
the statement of Theorem \ref{AH}.

At the end of \cite{Pa2} it is proved that the expression
of the general element of $\sigma_7(V^{5,2})$ has a sum
of seven $5$-th powers is unique. This fact was proved also by Richmond
$\cite{Ri}$, and also Hilbert knew and claimed it in a letter to Hermite in 1888 \cite{Hil}.
For recent results about the uniqueness of canonical forms see \cite{Me}.

The work of Terracini is a turning point in this story.
In his celebrated paper \cite{Te1} Terracini
introduced new techniques to attack the problem, and in particular
he proved (what today are called) the first and the second Terracini lemmas, as we have stated in Section \ref{notation}.
These results are not difficult to prove, but they represent a new viewpoint on the subject.
Terracini got them in an elegant way, as a natural state of things.
In \cite{Te1} Terracini was actually interested in a different direction.
Before of his work there were two different
characterizations of the Veronese surface.
Del Pezzo proved in 1887 that the Veronese surface in $\P 5$ is the unique surface
such that any two of its tangent planes
meet  each other. Severi proved in 1901  that the Veronese surface in $\P 5$
is the unique surface such that its
secant variety does not fill the ambient space. Is this only a coincidence?
Terracini's approach allows to unify these two results, indeed
thanks to the first Terracini lemma  the results of Del Pezzo and Severi turn out to be equivalent.
 This was probably not a surprise because the Severi proof
was deeply inspired by the Del Pezzo proof. But this opens another story that we do not pursue here.

In 1915 Terracini, with the paper \cite{Te2} realized that his two lemmas
allow to attack the problem raised by Palatini.
Terracini obtained in few lines at page 93 Theorem
\ref{AH} in the case $n=2$. His argument is the following.
The general ternary form of degree $d$ is sum of the expected number $k=\lceil \frac{(d+2)(d+1)}{6}\rceil $
of $d$-th powers of linear form if and only if there is no plane curve having
double points at general $p_1, \ldots p_k$.
On the other hand if there is such a curve, by
Lemma \ref{ter2} it has to contain as a component a double curve of
degree $2l$ through $p_1, \ldots p_k$. Hence we have the inequality
$$ k\le\frac{l(l+3)}{2}$$
so that we get the inequality
$$\left\lceil\frac{(d+2)(d+1)}{6}\right\rceil \le \frac{d}{4}({\frac{d}{2}}+3)$$
which gives $d=2$ or $d=4$ as we wanted.
This is the third published proof of Theorem \ref{AH} in the case $n=2$,
and the reader will notice that it is a refinement of Campbell proof.

Terracini observed also in \cite{Te2} that the exceptional case of cubics in $\P 4$ is solved by
the consideration that given seven points in $\P 4$, the rational quartic through them is the singular locus
of its secant variety, which is the cubic hypersurface defined by the invariant $J$ in the theory of binary quartics.

In \cite{Te3} Terracini got a proof of \thmref{AH} for $n=3$. In the
introduction he finally quoted the paper of Campbell, so it is
almost certain that he was not aware of it when he wrote the
article \cite{Te2}. Terracini gave to Campbell the credit to have
stated correctly Theorem \ref{AH} in the cases $n=2$ and $n=3$. We
quote from \cite{Te3}: {\it ``Questa proposizione fu dimostrata per la
prima volta in modo completo dal Palatini \cite{Pa2}, vedi
un'altra dimostrazione nella mia nota \cite{Te2}; ma gi\`a l'aveva
enunciata parecchi anni prima J.E.\ Campbell \cite{Ca} deducendola
con considerazioni poco rigorose, considerazioni che divengono
anche meno soddisfacenti quando il Campbell passa ad estendere la
sua ricerca alle forme quaternarie.''}
\footnote{\it This proposition was completely proved for the first time by
Palatini \cite{Pa2}, see another proof in my note \cite{Te2}; however
J.E.\ Campbell \cite{Ca} already stated it several years before,
deducing it in a not very rigorous way,
and his argument becomes even less satisfactory when Campbell
tries to extend his research to quaternary forms.}

This claim about the lack of rigor is interesting, because after a few years the Italian
school of algebraic geometry received the same kind of criticism,
especially from the Bourbaki circle.
The concept of the measure of rigor,  invoked by Terracini, is also interesting. Indeed we can
agree even today that Campbell argument was essentially correct
in the case $n=2$, but it was wrong in the case $n=3$.

Terracini's paper \cite{Te3} represents a change in the writing style.
All the lemmas and the theorems are ordered and numbered,
differently from all the papers quoted above.
His proof is by induction on the degree, and he uses
what we called in Section \ref{sezione4} the Castelnuovo sequence,
by specializing as many points as possible on a plane.
We saw in Section \ref{sezione4} that there is an arithmetic problem which
makes the argument hard when the number of double points is near to a critical bound.
Terracini's argument plays with linear systems with vanishing jacobian.
His approach was reviewed and clarified by Ro\'e, Zappal\`a and Baggio
in \cite{RZB}, during the 2001 Pragmatic School directed by Ciliberto and Miranda.
It seems to us that they also filled a small gap at the end of
Terracini's proof, obtaining a rigorous proof of
\thmref{AH} in the case $n=3$. It seems also that this approach does not
generalize to higher values of $n$.

In 1931 it appeared the paper \cite{Br} of Bronowski, at that time in Cambridge.
He took the statement of Theorem \ref{AH} from
\cite{Pa2} and he claimed to give a complete proof of it.
The argument of Bronowski is based on the possibility to check if a linear system
has vanishing jacobian by a numerical criterion.
This criterion already fails in the exceptional case
of cubics in $\P 4$, and Bronowski tried to justify this fact arguing that
the cases $n=2$ and $n=3$ are special ones.
However it is hard to justify his approach of considering
the base curve of the system.
In his nice MacTutor biography on the web, accounting a very active life, it is written:
{\it ``In 1933 he (Bronowski) published a solution of the classical functional
Waring problem, to determine the minimal $n$ such that a general degree $d$ polynomial $f$
 can be expressed as a sum of $d$-th powers of $n$
 linear forms, but his argument was incomplete.''}
We agree with this opinion.

In 1985 Hirschowitz \cite{Hir} gave a proof of Theorem \ref{AH}
 in the cases $n=2$ and $n=3$, which makes a step beyond the classical proofs, apparently not known to him at that time.
He used the powerful language of zero dimensional schemes in the degeneration argument,
this is the last crucial key to solve the general problem.
 In 1988 Alexander used the new tools introduced by Hirschowitz and in \cite{A} he  proved  \thmref{AH} for $d\ge 5$
with a very complicated but successful inductive procedure. He needed only a limited number
of cases for $d\le 4$ in the starting point of the induction.
In the following years Alexander and  Hirschowitz got \thmref{AH} for $d=4$ (\cite{AH0}) and finally
in \cite{AH1} they settled the case $d=3$, so obtaining the first complete proof of Theorem \ref{AH}.
  This proof,
which in its first version covered more than 150 pages, can be celebrated as a success of modern
cohomological theories facing with a long standing classical problem.   In 1993
Ehrenborg and Rota \cite{ER}, not aware of the work by Alexander and Hirschowitz, posed the problem
of Theorem \ref{AH} as an outstanding one.

In 1997 Alexander and Hirschowitz themselves got a strong simplification of
their proof in \cite{AH2}, working for $d\ge 5$.  By reading \cite{AH2}
it is very clear the role of the {\it dime} and the {\it degue}, see the Remark at the end of Section 6.
Later K.\ Chandler (see \cite{Ch}) simplified further the proof by Alexander
and Hirschowitz in the case $d\ge 4$, with the help of the Curvilinear \lemref{curvilinear}.
 In \cite{Ch2} she got a simpler proof also in
the case $d=3$.

Recently a different combinatorial approach to the problem
succeeded in the case $n=2$. The idea is to degenerate the
Veronese surface to a union of $d^2$ planes, as we learned from
two different talks in 2006 by R.\ Miranda and S.\ Sullivant. If
in the union of planes we can locate $k$ points on $k$ different
planes in such a way that the  corresponding planes are
transverse, then by semicontinuity the dimension of the $k$-secant
variety is the expected one. A proof of Theorem \ref{AH} in the
case $n=2$ along these lines was published by Draisma \cite{Dr}.
The proof reduces to a clever tiling of a triangular region. This
proof was extended to $n=3$ in S.\ Brannetti's thesis \cite{Bra}.
At present it is not clear if this approach, which is related to
tropical geometry, can be extended to $n\ge 4$.

We believe that the work on this beautiful subject will continue in the future.
Besides the higher multiplicity case
mentioned in the introduction, we stress that the equations of the higher secant
varieties $\sigma_k(V^{d,n})$ are still not known in general for $n\ge 2$,
and their knowledge could be useful in the applications.

\bigskip
{\sf
Maria Chiara Brambilla:\\
Dipartimento di Matematica e Applicazioni per l'Architettura, Universit\`a di Firenze\\
piazza Ghiberti 27,  50122 Firenze, Italy\\
brambilla@math.unifi.it\\

Giorgio Ottaviani:\\
Dipartimento di Matematica U. Dini, Universit\`a di Firenze\\
viale Morgagni 67/A,  50134 Firenze, Italy\\
ottavian@math.unifi.it
}
\end{document}